\documentclass[12pt]{article}
\usepackage{amsmath,amssymb,amsthm,latexsym}

\newtheorem{df}{Definition}
\newtheorem{lemma}{Lemma}
\newtheorem{corollary}{Corollary}
\newtheorem{theorem}{Theorem}
\newtheorem{prop}{Proposition}

\numberwithin{equation}{section}

\begin{document}
\title{{\sc Wrapping Brownian motion and heat kernels  II: symmetric spaces}}

\author{David G. Maher}
\maketitle

\newcommand{\g}{\mathfrak{g}}
\newcommand{\h}{\mathfrak{h}}
\newcommand{\s}{\mathfrak{s}}
\newcommand{\kg}{\mathfrak{k}}
\newcommand{\tg}{\mathfrak{t}}
\newcommand{\ug}{\mathfrak{u}}
\newcommand{\p}{\mathfrak{p}}
\newcommand{\X}{\mathfrak{X}}
\newcommand{\af}{\mathfrak{a}}
\newcommand{\D}{\mathbf{D}}
\newcommand{\R}{\mathbb{R}}
\newcommand{\C}{\mathbb{C}}
\newcommand{\Z}{\mathbb{Z}}
\newcommand{\N}{\mathbb{N}}
\newcommand{\E}{\mathbb{E}}
\newcommand{\bbP}{\mathbb{P}}
\newcommand{\T}{\mathbb{T}}
\newcommand{\F}{\mathcal{F}}
\newcommand{\A}{\mathcal{A}}
\newcommand{\x}{\mathbf{x}}
\newcommand{\y}{\mathbf{y}}
\newcommand{\ad}{\mathrm{ad}}
\newcommand{\Ad}{\mathrm{Ad}}
\newcommand{\Exp}{\mathrm{Exp}}
\newcommand{\grad}{\mathrm{grad} \,}
\newcommand{\hf}{\mathfrak{h}}
\newcommand{\cosec}{\mathrm{cosec} \,}
\newcommand{\csch}{\mathrm{csch} \,}
\newcommand{\sgn}{\mathrm{sgn} \,}
\newcommand{\hp}{\mathfrak{h_p}}

\begin{abstract}

In this paper we extend our previous results on wrapping Brownian
motion and heat kernels onto compact Lie groups to various symmetric
spaces, where a global generalisation of Rouvi\`ere's formula and
the $e$-function are considered. Additionally, we extend some of our
results to complex Lie groups, and certain non-compact symmetric
spaces.

\end{abstract}

{\small {\it Keywords:} Symmetric spaces, complex Lie group,
Brownian motion, wrapping map.

{\it AMS 2010 Subject classification:} 43A80, 22E30, 58J65.}

\section{Introduction}

In our previous paper \cite{M2}, we wrapped Brownian motion and heat
kernels from a compact Lie algebra (viewed as a Euclidean vector
space) to a compact Lie group using the wrapping map, $\Phi$, of
Dooley and Wildberger \cite{DW2}.  Recall that $\Phi$ was defined
for a suitable distribution, $\nu$, by
\begin{equation}
\langle \Phi(\nu), f \rangle = \langle \nu,j \tilde{f} \rangle
\end{equation}
where $f \in C^\infty (G)$, $\tilde{f} = f \circ \exp$ and $j$ the
analytic square root of the determinant of the exponential map. The
principal result is the {\bf wrapping formula}, given by
\begin{theorem}\label{thm2} (\cite{DW2}, Thm. 2) {Let $\mu, \nu$ be $G$-invariant distributions of compact support on $\g$ or two $G$-invariant integrable functions, then
\begin{equation}\label{wrapform}
\Phi(\mu * \nu) = \Phi(\mu) * \Phi(\nu)
\end{equation}
where the convolutions are in $\g$ and $G$, respectively.}
\end{theorem}

In this paper we consider wrapping Brownian motion in the context of
various symmetric spaces, where a global generalisation of the
wrapping formula \ref{wrapform} utilising Rouvi\`ere's $e$-function
is considered.\\

We will firstly recall the theory of $e$-functions in section 5,
then combine this with the theory of wrapping in section 6 where
global properties of the $e$-function are considered.  We then show
how to wrap heat kernels onto compact symmetric spaces from their
(Euclidean) tangent space, but due to the appearance of the
$e$-function, this requires some different ideas to those in our
previous work. Although the complicated nature of the $e$-function
makes explicit calculations difficult, our analysis provides several
insights into both compact and non-compact symmetric spaces.\\

Results concerning Brownian motion and heat kernels on symmetric
spaces have been previously given by many authors.  Our method
differs by using the wrapping map, which can be viewed as a global
version of the exponential map.   Thus, our results presented in
sections 7 and 8 are obtained are in the spirit of the tangent
space analysis advocated by Helgason (\cite{HE2}, \cite{HE3}).\\

In section 7 we consider the case of compact symmetric spaces, where
our work explains why the well-known Gaussian approximation  - also
known as the ``sum over classical paths" (see \cite{CAM},
\cite{DOW})  - does not give exact results for compact symmetric spaces that are not Lie groups.\\

In section 8 we consider the case of non-compact symmetric spaces,
we are able to quickly obtain the heat kernels for complex Lie
groups by wrapping.  We then discuss extensions of these results to
the non-compact symmetric spaces of ``split rank" type, in
particular the spaces $G/K$, $G$ complex, where our work implies
that the local convolution formula obtained by Torossian (\cite{T})
hold globally.

\section{Acknowledgements}

These results were obtained during the authors Ph.D candidature at
the University of New South Wales (see \cite{M}).  I would like to
thank my supervisor, Tony Dooley, for all his support and guidance
during this time.

\section{Notation and Formulae for Riemanniam Symmetric Spaces}

Let $G$ be a connected, semisimple Lie group, and $\g$ its Lie
algebra, with exponential map $\exp : \g \rightarrow G$.  Let
$\sigma$ denote an involutive automorphism of $G$ with respect to a
fixed compact subgroup $K$, such that we can write $\g = \kg \oplus
\p$ where $\kg$ and $\p$ are the eigenspaces of $\sigma$.\\

We identify $\p$ as the tangent space of the Riemannian symmetric
space $G/K$.  This is, furthermore, and exponential map $\Exp : \p
\rightarrow G/K$.  If $\pi$ is the canonical mapping $\pi : G
\rightarrow G/K$, then $\Exp = \pi \circ \exp$ . Let $\af$ denote a
maximal abelian subalgebra of $\p$, and call the dimension of $\af$
the {\it rank} of $G/K$.\\

Let $\Sigma = \Sigma(\g,\af)$ be the set of roots on $\g$ with
respect to $\af$, and $W(\g,\af)$ the corresponding Weyl group.  The
{\it restricted roots} of $\p$, denoted by $\Sigma_r$ are the roots
of $\g$ restricted to $\af$, and the {\it restricted Weyl group} is
the group generated by the reflections of the restricted roots.\\

Further details of restricted roots may be found in \cite{KNP} Ch.
VI.  We will use Knapp's definitions and notations for roots,
weights, etc.  We denote by $\af^+$ the positive Weyl chamber with
respect to the set of positive restricted roots, and $\bar{\af}^+$
its closure.  We also denote by $m_\alpha$ the multiplicity of the
root $\alpha$.  A restricted root $\beta$ is said to be {\it
multipliable} if there exists another restricted root $\alpha$ such
that $\beta = k \alpha$ for some integer $k \geq 2$. We denote
the set of multipliable roots by $\Sigma_m$.\\

The classification of Riemannian symmetric spaces can be found in
\cite{HE1}, Ch. V.  We will not give details here, but will briefly
outline the structure and relationship of the {\it compact type},
and the {\it non-compact type}:  Let $\g_\C$ be a complex Lie
algebra and $\ug$ a {\it compact real form}. $\g_\C$ may be regarded
as a real Lie algebra $\g^R = \ug \oplus i\ug$ with twice the (real)
dimension of $\ug$. We write $G$ for the Lie group corresponding
$\g_\C$, and $U$ for the Lie group corresponding to $\ug$. Thus, if
$(G,K)$ is a Riemannian symmetric pair, then $(U,K)$ is called its
{\it compact dual pair}, and $U/K$ is a compact Riemannian symmetric
space.  For example, let $U = SU(2)$ with $\ug = \mathfrak{su}(2)$.
Then $\g_\C = \mathfrak{sl}(2,\C)$ and $G = SL(2,\C)$.\\

Let $A = \exp(\af)$, and $A^+ = \exp(\bar{\af}^+)$, and Let $M$
denote the centraliser of $A$ in $K$.  Every element in $G$ has a
decomposition as $k_1 a k_2$, with $k_1 , k_2 \in K$ and $a \in A$.
In this decomposition, $a$ is uniquely determined up to conjugation
by a member of the Weyl group (see \cite{KNP} Thm. 7.39).
Additionally, $G$ has a Cartan decomposition $KA^+K$.  We have the
following integral formula for symmetric spaces (see \cite{HE2} Ch.
I, Thm. 5.8):
\begin{equation}\label{5214}
\int_{G/K} f(x) dx = c \int_{K/M} \biggl( \int_{A^+} f (ka \cdot o)
\delta (a) da \biggr) dk_M, \phantom{abcde} f \in C_c(G/K),
\end{equation}
 where $c$ is a suitable constant, $dx$ is the
$G$-invariant measure on $G/K$, $dk_M$ is the $K$-invariant measure
on $K/M$, normalised with mass 1, and
$$
\delta (\exp H) = \prod_{\alpha \in \Sigma_r^+} (\sinh \alpha
(H))^{m_\alpha}, \phantom{abcde} H \in \af^+.
$$

We note that $\Ad (k) \p \subseteq \p$, and that $\Ad-K$-invariant
functions are determined by their values on $\af^+$ (see \cite{HE2}
Ch. I, Thm. 5.17):
$$
\int_\p f(X) dX = c \int_{K/M} \biggl( \int_{\af^+} f (\Ad (k)H)
\delta_0 (H) dH \biggr) dk_M, \phantom{abcde} f \in C_c(\p),
$$
where $c$ is a suitable constant, $dk_M$ is the $K$-invariant
measure on $K/M$, normalised with measure 1, and
$$
\delta_0 (H) = \prod_{\alpha \in \Sigma_r^+} \alpha (H)^{m_\alpha},
\phantom{abcde} H \in \af.
$$

We denote by $J$ the Jacobian of $\Exp$,
$$
\int_{G/K} f(x) dx = \int_\p f(\Exp X)J(X)dX, \phantom{abcde} f \in
C_c (G/K)
$$
with $J$ given by $J = \delta / \delta_0$. We also write $J(X) = j^2
(X)$, where $j$ is calculated as:
$$
j(H) = \prod_{\alpha \in \Sigma_r^+} \biggl( \frac{\sinh
\alpha(H)}{\alpha (H)} \biggr)^{m_\alpha /2}, \phantom{abcde} H \in
\af.
$$

However, we will require that $j$ be smooth and real valued, which
is clearly not the case globally for every Riemannian symmetric
space. $j$ is smooth and real valued for compact Lie groups and
complex Lie groups (since their roots have even multiplicities), but
this is not so in general, and some of our
results will only be valid within a fundamental domain of the exponential map.\\

Let $\{ X_i \}_{i = 1}^n$ be an orthonormal basis of $\g$.  We
recall from the (left) action of an element of the Universal
enveloping algebra, $\mathfrak{U}(\g)$, on $C^\infty (G)$:
\begin{equation}
(f X)(g) = f(X ; g) = \biggl( \frac{\partial}{\partial t}
\biggr)_{t=0} f(g \exp^{t X}), \phantom{abcde} X \in \g, \; g \in G.
\end{equation}

We write {\bf D}$(\p)$ and $\D(G/K)$ for the set of left
$K$-invariant differential operators on $\p$ and $G/K$,
respectively. We also write $L_\p$ and $L_{G/K}$ for the Laplacians
on $\p$ and $G/K$, respectively.  Let $\mathfrak{U}(\g)^\kg$ be the
centraliser of $\kg$ in $\mathfrak{U}(\g)$.  The following results
can be found in \cite{HE2}, Ch II.  $\D(G/K)$ can be described in
terms of $\mathfrak{U}(\g)$ by
$$
\D(G/K) \cong \mathfrak{U}(\g)^\kg / (\mathfrak{U}(\g)^\kg \cap
\mathfrak{U}(\g) \kg)
$$

and $\D (\p)$ is isomorphic to $\D(\R^n)$.  We let $S(\p)$ denote
the symmetric algebra of $\p$.  For $\D(G/K)$ we have that $\D
(G/K)$ is a polynomial algebra, generated by algebraically
independent generators $D_1, \dots ,D_l$, whose degrees $d_1, \dots
,d_l$, where $l = rank(X)$, and are canonically determined by $G$.\\

When $X$ is of rank 1 we have that $\D(G/K)$ consists of polynomials
on $L_{G/K}$.  We recall the Casimir element of $G$, $\Omega$, is
defined by $\Omega = \sum_{i=1}^n X_i^2$.  In particular, consider
the orthonormal bases $(Y_j)$ and $(Z_k)$ (with respect to $B$) of
$\p$ and $\kg$, respectively.  Then we have
$$
\Omega = \Omega_\p - \Omega_\kg
$$
where $\Omega_\p = \sum_j Y_j^2$ and $\Omega_\kg = \sum_k Z_k^2$,
and moreover:
\begin{equation}
L_{G/K} f(xK) = f(\Omega ; x) = f(\Omega_\p ; x)
\end{equation}

We will also denote the transpose of a differential operator, $D$,
on $G/K$ by $D^t$.  We say that a differential operator $D$ is
symmetric if it satisfies $D^t = D$, and note that $L_{G/K}$ is a
symmetric operator.  The Laplacians $L_{G/K}$ and $L_\p$ are related
as follows (\cite{HE2} Ch. II, Prop. 3.15): the image
$L_{G/K}^{\Exp^{-1}}$ of $L_S$ under $\Exp^{-1}$ is given by:
$$
L_{G/K}^{\Exp^{-1}} f = (j^{-1} L_\p \circ j) f - j^{-1} (L_\p j) f
$$
for each $K$-invariant $C^\infty$ function $f$ on $\p$.  It is
important to note as per our above remarks on the $j$ function, this
result holds in general only in a fundamental domain of Exp, but
globally in the case of a compact or complex Lie group (\cite{HE2}
Ch. II, $\S 3$).\\

We now recall from \cite{GAN2} some of the theory of spherical
functions on Riemannian symmetric spaces. Let $G$ be a locally
compact group and $K$ a compact subgroup.  Let $\phi$ be a complex
valued function on $G/K$ with $\phi (o) = 1$. $\phi$ is said to be a
{\bf $K$-spherical function} on $G/K$ if it satisfies:
\begin{list}{}{}
\item (i) $\phi (k g K) = \phi(gK) \; \; \forall k \in K, \, g \in G$,
\item (ii) $D \phi = \lambda_D \phi \; \; \text{for each} \; D \in \D(G/K)$, \; \text{where} $\lambda_D \in
\C$.
\end{list}

Property $(i)$ ensures that a spherical function is determined by
its values on $A^+$.  We will also refer to these functions as {\bf
bi-$K$-invariant} functions.  Since a symmetric pair $(G,K)$ is
always a Gelfand pair (\cite{GV} Ch. I, Cor. 1.5.6), we may
construct the elementary spherical functions of $G/K$, denoted by
$\phi_\pi$, based on the representations of $G$ by
$$
\phi_\pi (x) := \langle \pi (x) e_\pi , e_\pi \rangle
$$
where $\pi$ is an irreducible unitary representation of $G$ on a
Hilbert space $\mathcal{H}$, and $e_\pi \in \mathcal{H}$ is a
$K$-fixed unit vector.  In the case where $U/K$ is a Riemannian
symmetric space of the compact type. The elementary spherical
functions for $U/K$ may also be calculated from the characters of
$U$
\begin{equation}  \varphi (g) = \int_K \chi (g^{-1} k) dk
\end{equation}
where $\chi$ is the character of an irreducible spherical
representation of $U$.\\

We retain the usual notations of the function spaces.  In
particular, we will consider the case when these functions are also
spherical functions  -  for example, we let $L^p (K \backslash G /
K)$ be the space of bi-$K$-invariant spherical $L^p$ functions on
$G/K$, and $\mathcal{S}(K \backslash G / K)$ be the set of
bi-$K$-invariant spherical Schwartz functions on $G/K$.  For $f \in
L^1 (K \backslash G / K)$ we define the Fourier transform
$\hat{f}(\lambda)$ by
$$
\hat{f}(\lambda) = \int_G f(x) \varphi_\lambda (x^{-1}) dx,
\phantom{abcde} \lambda \in \af^*.
$$
where $\varphi_\lambda$ is an elementary spherical function
corresponding to the weight $\lambda$.  Harish-Chandra described the
following Fourier inversion formula (see \cite{GV}, Thm. 6.4.2.) for
bi-$K$-invariant Schwartz functions:
$$
f(x) = \frac{1}{|W|} \int_{\Lambda} \hat{f}(\lambda) \varphi_\lambda
(x) |c(\lambda)|^{-2} d\lambda, \phantom{abcde} \lambda \in \af^*.
$$
The function $c(\lambda)$ on $\Lambda$ is such that
$c(\lambda)^{-1}$ is a tempered distribution, $|c(w \lambda)| =
|c(\lambda)|$ for $w \in W$, and we note that $|c(\lambda)|^{-2}
d\lambda$ is the Plancherel measure. We also recall that the
convolution of two $K$-invariant distributions, $\mu$ and $\nu$, on
$G/K$ is defined by:
\begin{equation}
(\mu * \nu, f) = (\mu (xK) \otimes \nu (yK), f(xyK) ), \phantom{abc}
x,y \in G,
\end{equation}
for any test function $f$ on $G/K$.

\section{Heat kernels and Brownian motion}

Heat equations on Riemannian symmetric spaces have been studied by
many authors in a variety of ways.  We define the heat equation by:
\begin{equation}
\frac{\partial}{\partial t} u(x,t) = L_{G/K} u(x,t)
\end{equation}
with initial data $u(x,0) = f(x)$.  The solution on the Cauchy
problem is given by
\begin{equation}
u(x,t) = \int_{G/K} h_t (x,y) f(y) dy
\end{equation}
where $h_t$ is the heat kernel.  We summarise some key properties of
$h_t$:
\begin{theorem} (c.f \cite{BER}, \cite{BGM}) $h_t$ satisfies the following properties: for all $x, y \in G/K$,
\begin{list}{}{}
\item (1) $h_t(x,y) = h_t(y,x) > 0$,
\item (2) $h_t$ is the density of a probability measure, with $\lim_{t \rightarrow 0} h_t (x,y) = \delta_x (y)$,
\item (3) $(\tfrac{\partial}{\partial t} - L_{G/K} ) h_t = 0$,
\item (4) $h_{t+s} (x,y) = \int_{G/K} h_t (x,z) h_s (z,y) dz$.
\end{list}
Moreover, if $U/K$ is compact, then $h_t$ can be expressed in terms
of the eigenvalues and eigenfunctions of the Laplacian as
$$
h_t(x,y) = \sum_{\lambda \in \Lambda^+} e^{-(\|\lambda+ \rho\|^2 -
\| \rho \|^2) t} \varphi_\lambda (x) \varphi_\lambda (y).
$$
\end{theorem}

These properties of $h_t$ can be shown to hold on more general
manifolds  -  the reader is referred to the listed sources.  On
$G/K$, the $G$-invariance implies that:
\begin{theorem} (c.f \cite{BER}, \cite{BGM}) $h_t$ also satisfies the following properties: for all $x, y \in G/K$,
\begin{list}{}{}
\item (1) $h_t(xK,yK) = h_t(y^{-1}x)$ is a convolution kernel,
\item (2) $x \mapsto h_t (x)$ is $K$-invariant on $G/K$, and
thus determined by its restriction to the positive Weyl chamber.
\end{list}
\end{theorem}

\noindent {\bf Remark:} These formulas can be used to derive the
heat kernel on a compact Lie group:
$$
H(g,t) = \sum_{\lambda \in \Lambda^+} d_\lambda e^{-(\|\lambda+
\rho\|^2 - \| \rho \|^2) t} \chi_\lambda (g), \phantom{abcde} g \in
G, \, t \in \R^+,
$$
which follows since the characters are the eigenfunctions of the
Laplacian, with eigenvalue $\|\lambda+ \rho\|^2 - \| \rho \|^2$, and
$\chi_\lambda (e) = d_\lambda$.\\

We now define Brownian motion on Riemannian manifolds: Suppose $M$
is an $n$-dimensional Riemannian manifold and $X_1 ,\dots, X_m$ are
vector fields on $M$.  If $(B_t)_{t \geq 0}$ is an $m$-dimensional
Brownian motion on $\R^n$ and $p \in M$, then an $M$-valued
stochastic process $(\xi_t)_{t \geq 0}$ is said to be a solution of
\begin{equation}
d\xi_t = \sum_{i=1}^n X_i (\xi_t) \circ dB_t^{(i)}, \phantom{abcde}
\xi_0 = p
\end{equation}
if for each $f \in C^\infty (M)$ we have
\begin{equation}\label{BMRM}
f(\xi_t) = f(p) + \sum_{i=1}^n \int_0^t (X_i f)(\xi_s) \circ
dB_s^{(i)}
\end{equation}

The solution of (\ref{BMRM}) is a {\bf Brownian motion} on $M$,
starting at $p \in M$.\\

For further details of Brownian motion, Stratonovich and It\^o
integrals on Riemannian manifolds and Riemannian symmetric spaces,
the reader is referred to \cite{LIAO2}, Ch. 2.  Importantly, we note
from this source that if $(B_t)_{t \geq 0}$ is a Brownian motion on
$G/K$, then $L_{G/K}$ is the generator of $(B_t)_{t \geq 0}$, and
$(B_t)_{t \geq 0}$ satisfies
\begin{equation}
\E(f(B_t)) = \int_{G/K} h_t (x,y) f(y) dy
\end{equation}

\section{Rouvi\`ere's $e$-function}

We now define a version of the wrapping map for Riemannian symmetric
spaces. Let $G/K$ be a Riemannian symmetric space with tangent space
$\p$. Let $\nu$ be a distribution of compact support on $\g$, and $f
\in C^\infty_c (G)$. We define the {\bf wrapping map}, $\Phi$ on
$G/K$ by
\begin{equation}
( \Phi(\nu), f )_{G/K} = ( \nu,j \cdot f \circ \Exp )_\p
\end{equation}
where $j$ the analytic square root of the determinant of the
exponential map.  We call $\Phi(\nu)$ the {\bf wrap} of $\nu$.\\

The case of a compact Lie group is particularly nice since the
wrapping map is a homomorphism between the convolution algebras of
Ad-invariant Schwartz functions or distributions of compact support
on $\g$, and central measures or distributions on $G$, ie,
\begin{equation}
\Phi(\mu *_\g \nu) = \Phi(\mu) *_G \Phi(\nu)
\end{equation}

In the case of a compact symmetric space, the wrapping map is no
longer a homomorphism: the convolution on the tangent space becomes
``twisted'' by a function denoted by $e$:
\begin{theorem}\label{efunct} (\cite{D1})  Let $U/K$ be a compact symmetric space with tangent space $\p$, and $\mu$ and $\nu$ $K$-invariant Schwartz functions or distributions of compact support
on $\p$.  There is a function $e : \p \times \p \rightarrow \C$ such
that
\begin{equation}\label{wrape1}
\Phi(\mu) *_{U/K} \Phi(\nu) = \Phi(\mu *_{\p, e} \nu)
\end{equation}
where
\begin{equation}\label{wrape2}
(\mu *_{\p, e} \nu)(X) = \int_\p \mu(Y)\nu(X-Y)e(X,Y)dY
\end{equation}
\end{theorem}

We will call the expression $\mu *_{\p, e} \nu$ a {\bf twisted convolution}.\\

The $e$-function was introduced by Rouvi\`ere in \cite{R1}, where a
local version of Theorem (\ref{efunct}) was proved for general
symmetric spaces (see \cite{R1}, Prop. 4.1).  For the case of the
2-sphere, the wrapping map and Rouvi\`ere's formula were studied in
the thesis of Chung \cite{CHU}.  Efforts to prove a global version
of Rouvi\`ere's formula for symmetric spaces have continued in
\cite{D1} and \cite{D2}.\\

The $e$-function arises from the following:  The left hand side of
(\ref{wrape1}) can be written as
\begin{align*}
(\Phi(\mu) *_{U/K} \Phi(\nu) , f ) & = \int_{G} \int_{G} \Phi (\mu)(xK) \Phi (\nu)(yK) f(xyK) dx dy\\
& = \int_\p \int_\p \mu(X) \nu(Y) j(X) j(Y) f(\exp X \Exp Y) dX dY .
\end{align*}

It is possible to show using the Campbell-Baker-Hausdorff series for
$\Exp$ (see \cite{D1} or \cite{R1}), that there exists $h, k \in \h$
such that
$$
\exp X \Exp Y = \Exp (h. X + k.Y) .
$$

We thus make the change of variables $(h. X, k.Y) \mapsto (X, Y)$,
and let the Jacobian of this transformation be $\psi (X,Y)$.  Our
expression then becomes
$$
\int_\p \int_\p \mu(h^{-1} X) \nu(k^{-1}Y) j(h^{-1}X) j(k^{-1}Y)
 f(\Exp (X + Y)) \psi (X,Y) dX dY .
$$

Since $j, \mu$ and $\nu$ are all $H$ invariant, this becomes
$$
\int_\p \int_\p \mu(X) \nu(Y) \frac{j(X) j(Y)}{j(X+Y)} \psi (X,Y)
(j.f \circ \Exp) (X + Y) dX dY .
$$

Putting
$$
e(X,Y) = \frac{j(X) j(Y)}{j(X+Y)} \psi(X,Y)
$$
we have
$$
\int_\p \int_\p \mu(X) \nu(Y) e(X,Y) (j.f \circ \Exp) (X + Y) dX dY
$$
which is (\ref{wrape2}). $\phantom{abcde} \square$\\

We will now briefly show how a global $e$-function may be
constructed for the case of the two-sphere.  This construction, due
to Rouvi\`ere, can be found in \cite{WLD} and \cite{D1} (the latter
eludicates as to how this may then be extended to all compact
symmetric spaces (\cite{D2})).\\

The global $e$-function is a ratio $g/f$ that compares the
convolution structures of $K$-orbits of $S^2$, to $K$-orbits of $\p
\cong \R^2$  -  the $K$-orbits in this case being circles centred at
the origin. To calculate $f$, consider two circles centred at the
origin of radius $r_1$ and $r_2$ on $\R^2$. For notational
convenience, it is best to consider these circles on the complex
plane.  We will consider the point $r_1$ on the first circle (we
could take any point, but we could obtain the same result by
rotation), and centre the circle of radius $r_2$ here.\\

We pick a point on the repositioned circle of radius $r_2$. This
point can be represented from the above construction as $r_1 + r_2
e^{i\theta}$, or as $r^{i\psi}$.  Therefore,
$$
r_1 + r_2 e^{i\theta} = r^{i\psi}
$$

We now vary the point on the circle of radius $r_2$ (by varying
$\theta$), and calculate how $r$ varies.  It is not hard to show
that
$$
\frac{2r}{2r_1 r_2 \sin \theta} \frac{dr}{\pi} = \frac{d\theta}{\pi}
$$

The denominator on the left-hand side is the area of the triangle on
the complex plane with vertices $0$, $r_1$ and $r_2 e^{i\theta}$.
By Heron's formula, we have:
$$
2 r_1 r_2 \sin \theta = \biggl( \prod_{\pm} (r \pm r_1 \pm r_2)
\biggr)^{1/2}
$$
where the product is taken over all choices of $+$ and $-$.  Thus,
the convolution of two circles of radius $r_1$ and $r_2$ has density
$$
f_{r_1, r_2} (r) = \frac{2r}{ \prod_{\pm} (r \pm r_1 \pm r_2)^{1/2}}
\chi_{[|r_1 - r_2|,r_1 + r_2]} (r)
$$

A similar calculation can be done on the surface of the two-sphere,
yielding:
$$
g_{r_1, r_2} (r) = \frac{\sin r}{\pi \sin r_1 \sin r_2} \prod_{\pm}
2 \sin \tfrac 12 (r \pm r_1 \pm r_2)^{1/2} \chi_{[|r_1 - r_2|,r_1 +
r_2]} (r)
$$

The $e$-function for the two-sphere is given by
$$
(g/f) (r) = \frac{\sin r}{\pi \sin r_1 \sin r_2} \prod_{\pm} \frac{
2 \sin \tfrac 12 (r \pm r_1 \pm r_2)^{1/2}}{(r \pm r_1 \pm
r_2)^{1/2}} \chi_{[|r_1 - r_2|,r_1 + r_2]} (r)
$$

In \cite{D2}, it will be shown that a global version of the
$e$-function for all compact symmetric spaces exists.  The proof
consists of reducing the calculation to the two-dimensional case and
using the above ideas.  The $e$-function for the $n$-dimensional
sphere is:
\begin{equation}\label{ecss}
e (X,Y) = \frac{\sin r}{\pi \sin r_1 \sin r_2} \biggl( \prod_{\pm}
\frac{ 2 \sin \tfrac 12 (r \pm r_1 \pm r_2)^{1/2}}{(r \pm r_1 \pm
r_2)^{1/2}}  \chi_{[|r_1 - r_2|,r_1 + r_2]} (r) \biggr)^{(n-3)/2}
\end{equation}
where $X \in \p$ is conjugate to $r_1 H$ in $\af$, $Y \in \p$ is conjugate to $r_2 H$ in $\af$, and $X + Y \in \p$ is conjugate to $r H$ in $\af$.\\

Since the wrapping formula is now a ``twisted" homomorphism, it is
no longer clear that we may wrap Brownian motion and the heat kernel
without some modification.  In the next section, we will show how
the Laplacian and $e$-function interact.

\section{Rouvi\`ere's formulae and the wrapping map}

We saw in \cite{M2} that the wrap of the Laplacian determines how
Brownian motion wraps. However, for general symmetric spaces we do
not have such a straightforward expression as we did in \cite{M2},
due to the twisted convolution involving the $e$-function (\ref{wrape1}).\\

The relationship between the $e$-function and differential operators
is given by Rouvi\`ere in \cite{R2}.  Rouvi\`ere expresses the
$e$-function as an infinite series and shows that it converges
 within a certain neighbourhood of $0 \in \p$.  Therefore, the results in \cite{R2}
 concerning the relationship between the $e$-function and differential
 operators only hold within this neighbourhood.\\

Equipped with the results from \cite{D1} and \cite{D2} presented in
the previous section concerning the global existence of the
$e$-function, we now show that Rouvi\`ere's formulae hold at least within a
neighbourhood of $0 \in \p$ where the $j$ function is smooth and real valued.\\

We will retain the notations and conventions given by Rouvi\`ere in
$\cite{R1}$ and $\cite{R2}$: let $G/K$ be a symmetric space, with
tangent space $\p$ and $\Exp: \p \rightarrow G/K$ the exponential
map, and $\p'$ be an $K$-invariant neighbourhood of $0 \in \p$ such
that $\Exp$ is a diffeomorphism.  Let $\Omega_0$ be an open
neighbourhood of $(0,0) \in \p' \times \p'$ satisfying the condition
that if $(X,Y) \in \Omega_0$, then $(k \cdot X, k \cdot Y) \in
\Omega_0$ for all  $k \in K$, and $(-X,-Y) \in \Omega_0$.  Let
$\p''$ denote the set of $X \in \p'$ such that $(X,0) \in
\Omega_0$.\\

\noindent {\bf Remark:}  The one exception to our notation is that
Rouvi\`ere (\cite{R1}, \cite{R2}) uses the notation $S$ to denote a
symmetric space and $\s$ instead of $\p$ as his formulae hold for
symmetric spaces which need not be Riemannian. We will continue to
use $G/K$ and $\p$ as we will only be considering Riemannian symmetric
spaces.\\

Since the set of $K$-invariant differential operators on $\p$ and
$G/K$ respectively, denoted by $\D(\p)$ and $\D(G/K)$ resp., are
polynomial algebras generated, respectively, by $L_\p$ and
$L_{G/K}$, we follow Rouvi\`ere's notation, we write the elements of
$\D (\p)$ of polynomials on
$\p^*$ as $p(\partial_X)$, $p$, or $p(\xi)$.\\

We write $\partial_\xi$ for $\partial / \partial_\xi$.  Let  $e_p
(X) = (p(\partial_Y)e)(X,0)$, and $\Delta_Y e (X) = \Delta_Y
e(X,0)$, where $\Delta_Y$ is the Laplacian acting on the second
variable, evaluated at $0$.  Let $\delta_o$ and $\delta_0$ be the
Dirac deltas on $G/K$ and $\p$, respectively. The map $p \rightarrow
\tilde{p}$ from is a linear isomorphism $\D(\p) \rightarrow
\D(G/K)$, according to
$$
(p^t \delta_0)^{\tilde{}} = \tilde{p}^t \delta_o
$$
or, more generally,
$$
\alpha * (p^t \delta_0)^{\tilde{}} = \tilde{p}^t \alpha
$$
for any distribution $\alpha$ on $G/K$.\\

Rouvi\`ere uses the notation $\tilde{f}$ as the inverse of the map
$f \mapsto j (f \circ \exp)$.  For the wrapping map (restricted to
the fundamental domain) this would read:
$$
\langle \Phi (\mu) , \tilde{f} \rangle = \langle \mu , f \rangle
$$

Our main result provides an analogue of \cite{M} Proposition 1 for
symmetric spaces:
\begin{prop}\label{eandlap2} Suppose $\mu$ is a $K$-invariant Schwartz functions
on $\p$, then on a suitable neighbourhood of $0 \in \p$,
$$
\Phi \bigl( (L_\p - \Omega_*) \mu \bigr) = L_{G/K} \bigl( \Phi (\mu)
\bigr)
$$
\end{prop}

The quantity $\Omega_*$ will be defined below, as well as what
constitutes a ``suitable neighbourhood''.  We will prove Proposition
\ref{eandlap2} by restating a number of results concerning
differential operators from \cite{R2} in the language of wrapping.\\

Firstly, we recall Theorem \ref{efunct}: If $\mu$ and $\nu$ are
$K$-invariant Schwarz functions or distributions of compact support
on $\p$, then
\begin{equation}\label{615}
\Phi(\mu) *_{G/K} \Phi(\nu) = \Phi(\mu *_{\p, e} \nu)
\end{equation}
where
\begin{equation}\label{616}
(\mu *_{\p, e} \nu)(X) = \int_\p \mu(Y)\nu(X-Y)e(X,Y)dY
\end{equation}

(\ref{615}) and (\ref{616}) may be written as follows (see \cite{R2}
Thm. 2.1, though note we have no restrictions on the supports of
$\mu$ and $\nu$):
\begin{equation}\label{617}
\langle \Phi(\mu) *_{G/K} \Phi(\nu) , f \rangle = \langle \mu(Y)
\nu(X), e(X,Y) j(X + Y) f(\Exp (X + Y) \rangle
\end{equation}

(\ref{617}) may be applied to give a relationship between wrapping,
differential operators, and the $e$-function.  We recall the
definition of the {\it symbol} of a differential operator on $\p
\cong \R^n$, given by:
$$
p(X,\partial_X) f(X) = \int_{\R^n} \biggl( \int_{\R^n} p(X,\xi)
\hat{f}(\xi) e^{2\pi i \langle Y, \xi \rangle} dY \biggr) d\xi
$$

Rouvi\`ere calculates the symbol for general differential operators
using the $e$-function:
\begin{theorem}\label{eanddiff} (\cite{R2} Thm. 3.1)  Let $G/K$ be a symmetric space, and $p \in I(\p)$, and let $U$ be an $K$-invariant open subset of
$\p''$.\\

(i) For any $K$-invariant distribution $u$ on $U$, one has
$$
\tilde{p}^t \tilde{u} = (( p^t(X, \partial_X))u)^\sim
\phantom{abcde} on \; \Exp U
$$
where $p(X, \xi)$ is the differential operator with analytic
coefficients on $\p''$ corresponding to the symbol
$$
p(X, \partial_X) = e(X, \partial_\xi) p(\xi) = \sum \frac{1}{\alpha
!}\partial_Y^\alpha e(X,0). \partial^\alpha_\xi p(\xi);
$$
here, $X \in \p''$, $\xi \in \p^*$ and $\sum$ is a (finite) summation over $\alpha \in \N^n$.\\

(ii)  If $f$ is any $H$-invariant $C^\infty$ function on $U$, one
has
$$
\tilde{p} \tilde{f} = (( p(X, \partial_X))f)^\sim \phantom{abcde} on
\; \Exp U
$$
\end{theorem}

We will be specifically considering the Laplacian, where Rouvi\`ere
proves:
\begin{corollary}\label{eandlap} (\cite{R2}, Cor. 3.6) Suppose $G/K$ is a Riemannian symmetric space with tangent space $\p$, we have\\

(i)  $\tilde{L_\p} = L_{G/K} + L_\p j(0)$,\\

(ii)  $L_{G/K} \tilde{u} = (L_\p u - j^{-1}L_\p j \, . u)\tilde{}$ on $\Exp W$.\\

Here, $W$ is an $K$-invariant subset of $\p''$, and $u$ is any $K$-invariant function, or distribution, on $W$.\\

(iii)  $j^{-1} L_\p j (X) = L_\p j (0) - L_Y e(X,0)$ on $\p''$,\\

(iv)  If $G$ is complex semisimple, then $L_\p j(0) = n/12$, where
$n = \dim G/K$.
\end{corollary}

We now state our version of Theorem \ref{eanddiff} (with proof
almost identical to Rouvi\`ere's) using the notation of wrapping:
\begin{theorem}\label{thm614}  For any $K$-invariant Schwartz functions or distributions of compact support $\mu$ on $\p$, we
have for a suitable neighbourhood of $0 \in \p$:
\begin{equation}\label{618}
\Phi(p) \Phi(\mu) = \Phi( (p(X,\partial_X))^t \mu)
\end{equation}
where $p(X,\partial_X)$ is the differential operator with analytic
coefficients on $\p$ corresponding to the symbol
$$
p(X,\xi) = e(X,\partial_\xi)p(\xi) = \sum_\alpha \frac{1}{\alpha !}
\partial^\alpha_Y e(X,0) \cdot
\partial^\alpha_\xi p(\xi)
$$

Moreover, if $f \in C^\infty_K (G/K)$, then
$$
\Phi(p) \Phi(f) = \Phi( (p(X,\xi))^t f)
$$
\end{theorem}

The next result that we require is:
\begin{prop}\label{prop615} (\cite{R2} Prop. 3.5(i)) Suppose $q$ is a
second order homogeneous operator.  Then, with the above notation,
\begin{equation}\label{619}
p(X,\partial_X) = p(\partial_X) + e_p (X) = (p(X,\partial_X) )^t
\end{equation}
\end{prop}


We now state our version of Corollary \ref{eandlap}, and define the
``suitable neighbourhood'' for which it holds:
\begin{corollary} Suppose $S$ is a
Riemannian symmetric space with tangent space $\p$, and $u$ is any
$K$-invariant function, or distribution, on $\p$.  Then for a
suitable neighbourhood of $0 \in \p$ comprising of the domain where
$j$ is smooth and real valued, we have:
$$
L_{G/K} \Phi(u) = \Phi \bigl( (L_\p - j^{-1}L_\p j )(u) \bigr)
$$

Moreover,
$$
j^{-1} L_\p j (X) = L_\p j (0) - L_Y e(X,0)
$$
\end{corollary}

\noindent {\bf Proof:}  Firstly, recall that the Laplacian is a
symmetric operator.  By Theorem \ref{thm614} and Proposition
\ref{prop615} we have
$$
\Phi (L_\p \mu + (\Delta_Y e(X,0) - (L_\p j)(0)) \mu ) = L_{G/K}
\Phi (\mu)
$$
where $j$ is smooth and real valued.  Taking $\mu = j$, by
\cite{DW1} Prop. 2.4, we have $\Phi(j) = 1$, and so
$$
\Phi (L_\p j + (\Delta_Y e(X,0) - (L_\p j)(0)) j ) = L_{G/K} 1 = 0
$$
and thus
$$
j^{-1} L_\p j = (L_\p j)(0) - \Delta_Y e(X,0)
$$
which concludes the proof. $\phantom{abc} \square$\\

In the next section, we will compute the term $j^{-1} L_\p j$, which
in turn will shed more light on the domain where our results hold.

\section{The Compact Case}

We now present our results for compact symmetric spaces , that is,
we show how to wrap Brownian motion onto a compact symmetric space
$U/K$ from its tangent space $\p$ and thus compute the heat kernel.
However, the situation is more complicated since the wrapping map is
no longer a homomorphism. We shall analyze this situation using
Rouvi\`ere's $e$-functions.\\

We are able to make some general statements about wrapping heat
kernels to $U/K$, but we are unable to explicitly compute the heat
kernels due to complicated potential terms arising from the
$e$-function.  These complications enable us to show why the {\it
Gaussian approximation} to the heat kernel does not give exact
results for compact symmetric spaces that are not compact Lie
groups.

\subsection{The Gaussian approximation}\label{section62}

Analysis of the heat kernel $K_t$ on a manifold $M$ of dimension $n$
by considering the heat kernel on its tangent space, combined with
information about the exponential map, has been considered by many
authors. For small time values, the {\it Minakshisundaram-Pleijel
(M-P) expansion} provides an expansion for the heat kernel on a
small neighbourhood of $M$: if $x$ and $y$ are close, then as $t
\rightarrow 0$,
\begin{equation}
K_t (x,y) = (2\pi t)^{-n/2} \exp \biggl( \frac{-d(x,y)^2}{2t}
\biggr) \times \Bigl( c_0 + t c_1 + \dots t^n c_n + O(t^{n+1})
\Bigr)
\end{equation}
where $d(x,y)$ is the Riemannian distance, and where the $c_k$'s are
dependent on $x$ and $y$. Further details of this expansion and
related techniques can be found, for example, in \cite{CAM}.\\

These neighbourhoods are then patched together to approximate the
heat kernel on $M$. This involves using the first order
approximation M-P expansion - the term $c_0$ can be shown to be
equal to $j^{-1}$. The solution which this technique produces is
sometimes referred to as the {\it Gaussian approximation} to the
heat kernel on $M$. It appears to have been first written down in
the thesis of Low \cite{LOW} (see also Camporesi \cite{CAM}, $\S
5$).
\begin{df} The {\bf Gaussian approximation} to the heat kernel on
$M$ is given by
\begin{equation}\label{Gapprox}
K_{\textrm{Gaussian}} (\Exp X,t) = \sum_{\gamma \in \Gamma}
\frac{p_t}{j} (X + \gamma), \phantom{abcde} X \in M.
\end{equation}
where $p_t$ is the heat kernel on $\R^n$.
\end{df}

In Dowker \cite{DOW} (see also Camporesi \cite{CAM}) the exactness
of the Gaussian approximation (\ref{Gapprox}) (except for a
so-called ``phase factor'' of $e^{\| \rho \|^2 t}$) is asserted to
hold for compact Lie groups, but not for general compact symmetric
spaces.\\

In the case of compact Lie groups, we assert that the reason why the
Gaussian approximation (\ref{Gapprox}) gives an exact expression of
the shifted heat kernel, but not in the case of general compact
symmetric spaces, is explained by the wrapping map and the $e$-function.\\

Firstly, we compute the wrap of a $K$-invariant Schwartz function as
a sum over the integer lattice. The proof is almost identical to the
proof of \cite{DW2}, Thm. 1:
\begin{prop}\label{614}  Suppose $\mu$ is a $K$-invariant Schwartz function on $\p$, given on $\af$.  Then,
$$
\Phi(j \mu)(\exp H) = \sum_{\gamma \in \Gamma} \mu (H + \gamma),
\phantom{abcde} H \in \af.
$$
\end{prop}
{\bf Proof:}  Let $\Psi$ be the $K$-invariant $C^\infty$ function on
$U/K$ given on $A$ by $\Psi (\Exp H) = \sum_{\gamma \in \Gamma} \mu
(H + \gamma)$.  For $f \in C^\infty (U/K)$,
\begin{align*}
(\Phi (j \mu), f) & = (j \mu , j \tilde{f}) = \int_\p j^2 (X) \mu (X) \tilde{f}(X) dX\\
& = \int_{\af^+} \prod_{\alpha \in \Sigma_r^+} \alpha^{m_\alpha} (H) j^2 (H) \mu (H) \int_{K} \tilde{f}(h \cdot H) dh dH\\
& = \int_{\af^+} \prod_{\alpha \in \Sigma_r^+} |\sin \alpha (\exp H)|^{m_\alpha} \mu (H) (\tilde{f})^{K}(H) dH\\
& = \frac{1}{|W|} \int_{\af} \prod_{\alpha \in \Sigma_r^+} |\sin
\alpha (\Exp H)|^{m_\alpha} \mu (H) (f^{K})^\sim (H) dH.
\end{align*}

If $\af_\Gamma \subseteq \af$ is a fundamental domain for $\Gamma$
in $\af$, then this becomes
\begin{align*}
\frac{1}{|W|} & \int_{\af_\Gamma} \prod_{\alpha \in \Sigma_r^+} |\sin \alpha (\Exp H)|^{m_\alpha} (f^{K})^\sim (H) \sum_{\gamma \in \Gamma} \mu (H + \gamma) dH\\
& = \frac{1}{|W|} \int_{A} \prod_{\alpha \in \Sigma_r^+} |\sin \alpha (a)|^{m_\alpha} (f^{K}) (a) da\\
& = \int_{U/K} \Psi(h) f(h) dh
\end{align*}
as required.  $\phantom{abc} \square$\\

Wrapping the heat kernel $p_t$ on $\p \cong \R^n$ yields the
Gaussian approximation (\ref{Gapprox}).  From \cite{M2}, in the case
of a compact Lie group we have
$$
\Phi (p_t) = q_t^\rho.
$$
In this case, the Gaussian approximation is the shifted heat kernel.
Moreover, note that the wrapping formula gives:
$$
\Phi (p_{t+s}) = \Phi (p_t * p_s) = \Phi (p_t) * \Phi(p_s) =
q_t^\rho * q_s^\rho = q^\rho_{t + s}.
$$

However, for compact symmetric spaces that are not Lie groups, we
have a non-trivial $e$-function that is ``twisting" the convolution
structure on $\p$, which ``twists'' wrapping the heat convolution
semigroup:
\begin{equation}\label{twistheat}
q_{t+s}^\rho = q_t^\rho *_{U/K} q_s^\rho = \Phi (p_t) *_{U/K} \Phi
(p_s) = \Phi (p_t *_{\p,e} p_s)
\end{equation}
but $\Phi (p_t *_{\p,e} p_s)$ is not equal to $\Phi (p_{t+s})$ when
$e$ is not identically 1.  Thus, the underlying reason that the
Gaussian approximation is not exact in this case is that we have a
twisted convolution on $\p$.  In summary we have:
\begin{theorem} The Gaussian approximation (\ref{Gapprox}) is not
equal to the heat kernel (modulo the phase factor of $e^{\| \rho
\|^2 t}$) for compact symmetric spaces that have a non-trivial
$e$-function.
\end{theorem}

Although we have a twisted convolution on the tangent space,
Rouvi\`ere's results from section 7 show that the potential term and
the $e$-function are closely related. That is,
$$
\Phi \bigl( (L_\p - \| \rho \|^2 + L_Y e(X,0) ) \mu \bigr) = L_{U/K}
\bigl( \Phi (\mu) \bigr)
$$
and
$$
- \| \rho \|^2 + L_Y e(X,0) = (j^{-1} L_\p j)(X)
$$

We now extend the results from \cite{M2} to wrapping Brownian motion
and the heat kernel onto $U/K$ from $\p$.  This involves a more
complex potential term than the constant $\| \rho \|^2$ term we
encountered in our previous results.\\

We now wrap Brownian motion in the same way that we did in Section
3.3, by putting $ h = j.f \circ \exp \in C^\infty_c (\p)$ in the
It\^o equation.  That is
$$
h(\zeta_t) = h(0) + \sum_{i=1}^n \int_0^t \sum_{i=1}^n
\frac{\partial h}{\partial x_i} (\zeta_t) dB_t^{(i)} + \tfrac 12
\int_0^t \bigl( L_\p + (j^{-1} L_\p  j) \bigr) (\zeta_t) dt
$$
wraps to
$$
f(\xi_t) = f(e) + \sum_{i=1}^n \int_0^t (X_i f)(\xi_t) dB_t^{(i)} +
\tfrac 12 \sum_{i=1}^n \int_0^t (L_{U/K} f)(\xi_t)dt
$$

Therefore, in light of this and Corollary \ref{eandlap2} we need to
consider the heat equation with potential $j^{-1} L_\p j$ given on
$\p \cong \R^n$ by
\begin{equation}\label{peturbHK}
\frac{\partial p_t^e}{\partial t} = (L_\p + (j^{-1} L_\p  j)) \,
p_t^e
\end{equation}
and we shall refer to the fundamental solution of (\ref{peturbHK}),
$p_t^e$, as the {\it perturbed heat kernel}. Wrapping the perturbed
heat kernel from $\p$ will then yield the heat kernel on $U/K$ since
these correspond to Laplacians on $\p$ and $U/K$:
\begin{theorem}\label{peturbHK2}  The wrap of the perturbed heat kernel $p_t^e$ on $\p$, given by
$$
\Phi(p_t^e)(\Exp H) = \sum_{\gamma \in \Gamma} \frac{p_t^e}{j} (H +
\gamma)
$$
is the the heat kernel on $U/K$.
\end{theorem}

However, the potential term $j^{-1} L_\p j$ is complex - even for
the case of the two-sphere, and finding the perturbed heat kernel
appears extremely difficult. As a first step, we now calculate these
potentials.

\subsection{The quantity $j^{-1} L_\p j$ for the compact case}

In this section we calculate $j^{-1} L_\p  j$.  In \cite{HE3} this
is referred to as the quantity $\Omega_*$ (but is not explicitly
calculated).  This calculation is similar in spirit to that of
\cite{HS}, section2.  In this paper, the authors calculate $(j \circ
\log) L_{G/K} (j \circ \log)^{-1}$ for the non-compact case.
However, this quantity in the compact case is not straightforward to
handle, as the singular values of $j^{-1}$ would require careful
treatment of $(j \circ \log) L_G (j \circ \log)^{-1}$, and certainly
would not lend itself easily to global analysis.\\

Firstly, note that we are simply taking the derivatives of a
function that is invariant under the action of $K$. Thus, to compute
$L_\p  j$, we firstly recall the following result from Helgason. Let
$\{ H_i \; : \; i = 1, \dots , l \}$ denote an orthnormal basis of
$\af$, and $D_{H_i}$ the first-order differential operator in the
$H_i$ direction, with $L_\af$ is the Laplacian on $\af$.  We have:
\begin{prop} (\cite{HE2}, Ch. 2, Prop 3.13)  For the adjoint action of $K$ on
$\p$ with transversal manifold $\af^+$, the radial part of the
Laplacian $L_\p$ is given by
\begin{align*}
\Delta(L_\p) & =\sum_{i=1}^l \biggl( D^2_{H_i} + \sum_{\alpha \in
\Sigma_r^+} m_\alpha \frac{1}{\alpha (\cdot)}
\alpha (H_i) D_{H_i} \biggr)\\
& = L_\af + \sum_{\alpha \in \Sigma_r^+} m_\alpha \frac{1}{\alpha}
A_\alpha
\end{align*}
\end{prop}

Here, the vector $A_\alpha$ is determined by $\langle A_\alpha , H
\rangle = \alpha (H), \; H \in \af$, and is considered as a
differential operator on $A^+ \cdot o$.
\begin{lemma}\label{632}  (\cite{OP} Thm. 5.1)  On the open dense subset of $\af$ where they are defined,
\begin{list}{}{}
\item (i) $$ \sum_{\alpha \neq k \beta, \; \alpha , \beta \in \Sigma_r^+} m_\alpha m_\beta \langle \alpha , \beta \rangle \frac{1}{\alpha \beta} =
0$$
\item (ii) \; \; In the case of a noncompact symmetric space $G/K$ we have:
$$ \sum_{\alpha \neq k \beta, \; \alpha , \beta \in \Sigma_r^+}
m_\alpha m_\beta \langle \alpha , \beta \rangle (\coth \alpha \coth
\beta - 1) = 0 $$
\item (iii) \; \; In the case of a compact symmetric space $U/K$ we have:
$$ \sum_{\alpha \neq k \beta, \; \alpha , \beta \in \Sigma_r^+}
m_\alpha m_\beta \langle \alpha , \beta \rangle (\cot \alpha \cot
\beta + 1) = 0 $$
\end{list}
\end{lemma}

We now use these expressions to calculate the quantity $(j(H)^{-1}
\Delta (L_{\p}) j)(H)$, which we will compute for both compact and
non-compact symmetric spaces.
\begin{prop}  With the above notation,
$$
(j^{-1} L_\p j) (H) = -\| \rho \|^2 + F(H)
$$
where
\begin{align*}
F(H) = 2\sum_{\alpha \in \Sigma_m^+} \frac{m_\alpha m_{2\alpha}}{4} \biggl( \cot \alpha (H) \cot 2\alpha (H) \biggr) | \alpha |^2 \\
 + \sum_{\alpha \in \Sigma_r^+} \frac{m_\alpha(m_\alpha - 2)}{4} \biggl( \cosec^2 \alpha (H) - \frac{1}{\alpha (H)^2} \biggr) | \alpha |^2
\end{align*}
\end{prop}

{\bf Proof:} Direct computation of the partial derivatives of $j$
yields:
$$
D_{H_i} j(H) = j(H) \sum_{\alpha \in \Sigma_r^+} \frac{m_\alpha}{2}
\biggl( \cot \alpha (H) - \frac{1}{\alpha(H)} \biggr) \alpha(H_i)
$$
\begin{align*}
\sum_{i=1}^l D^2_{H_i} j(H) = j(H) \Biggl[ \sum_{\alpha, \, \beta \in \Sigma_r^+} \frac{m_\alpha m_\beta}{4} & \biggl( \cot \alpha (H) - \frac{1}{\alpha(H)} \biggr) \biggl( \cot \beta (H) - \frac{1}{\beta (H)} \biggr) \langle \alpha, \beta \rangle \\
& - \sum_{\alpha \in \Sigma_r^+} \frac{m_\alpha}{2} \biggl( \cosec^2
\alpha (H) - \frac{1}{\alpha(H)^2} \biggr) | \alpha |^2 \Biggr]
\end{align*}

$$
\sum_{i=1}^l \sum_{\alpha \in \Sigma_r^+} m_\alpha \frac{1}{\alpha
(H)} \alpha (H_i) D_{H_i} j(H) = j(H) \sum_{\alpha, \, \beta \in
\Sigma_r^+} \frac{m_\alpha m_\beta}{2} \frac{1}{\alpha (H)} \biggl(
\cot \beta (H) - \frac{1}{\beta (H)} \biggr) \langle \alpha, \beta
\rangle
$$

Multiplying out and collecting like terms gives
\begin{align*}
\Delta (L_{\p}) j(H) = j(H) \Biggl[ \sum_{\alpha, \, \beta \in \Sigma_r^+} \frac{m_\alpha m_\beta}{4} & \biggl( \cot \alpha (H) + \frac{1}{\alpha(H)} \biggr) \biggl( \cot \beta (H) - \frac{1}{\beta (H)} \biggr) \langle \alpha, \beta \rangle \\
& - \sum_{\alpha \in \Sigma_r^+} \frac{m_\alpha}{2} \biggl( \cosec^2
\alpha (H) - \frac{1}{\alpha(H)^2} \biggr) | \alpha |^2 \Biggr]
\end{align*}

That is,
\begin{align*}
\Omega_* (H) & = (j^{-1} \Delta (L_{\p}) j)(H)\\
& = \sum_{\alpha, \, \beta \in \Sigma_r^+} \frac{m_\alpha m_\beta}{4} \biggl( \cot \alpha (H) + \frac{1}{\alpha(H)} \biggr) \biggl( \cot \beta (H) - \frac{1}{\beta (H)} \biggr) \langle \alpha, \beta \rangle \\
& \phantom{aaaaa} - \sum_{\alpha \in \Sigma_r^+} \frac{m_\alpha}{2} \biggl( \cosec^2 \alpha (H) - \frac{1}{\alpha(H)^2} \biggr) | \alpha |^2\\
& = \sum_{\alpha, \, \beta \in \Sigma_r^+} \frac{m_\alpha m_\beta}{4} \biggl( \cot \alpha (H) \cot \beta (H) + \cot \beta (H) \frac{1}{\alpha(H)}\\
& \phantom{aaaaabbbbbcccccdddddeeeeefffffggggg} - \cot \alpha (H) \frac{1}{\beta (H)} - \frac{1}{\alpha (H)} \frac{1}{\beta (H)} \biggr) \langle \alpha, \beta \rangle \\
& \phantom{aaaaa} - \sum_{\alpha \in \Sigma_r^+} \frac{m_\alpha}{2}
\biggl( \cosec^2 \alpha (H) - \frac{1}{\alpha(H)^2} \biggr) | \alpha
|^2
\end{align*}

Decomposing the first sum into diagonal and off-diagonal parts we
have
\begin{align*}
\Omega_* (H) & = \sum_{\alpha, \, \beta \in \Sigma_r^+, \alpha \neq \beta} \frac{m_\alpha m_\beta}{4} \biggl( \cot \alpha (H) \cot \beta (H) + \cot \beta (H) \frac{1}{\alpha(H)}\\
& \phantom{aaaaabbbbbcccccdddddeeeeefffffggggg} - \cot \alpha (H) \frac{1}{\beta (H)} - \frac{1}{\alpha (H)} \frac{1}{\beta (H)} \biggr) \langle \alpha, \beta \rangle \\
& \phantom{aaaaa} + \sum_{\alpha \in \Sigma_r^+}
\frac{m_\alpha^2}{4} \biggl( \text{cosec}^2 \alpha (H) - 1 -
\frac{1}{\alpha (H)^2} \biggr) | \alpha |^2
\end{align*}
\begin{equation}\label{sphere631}
\phantom{aaaaabbbbbcccccdddddeeeee} - \sum_{\alpha \in \Sigma_r^+}
\frac{m_\alpha}{2} \biggl( \text{cosec}^2 \alpha (H) -
\frac{1}{\alpha(H)^2} \biggr) | \alpha |^2
\end{equation}

Using Lemma \ref{632} (i) on the first sum, and combining the last
two, we have
\begin{align*}
\Omega_* (H) & = \sum_{\alpha, \, \beta \in \Sigma_r^+, \alpha \neq \beta} \frac{m_\alpha m_\beta}{4} \biggl( \cot \alpha (H) \cot \beta (H) + \cot \beta (H) \frac{1}{\alpha(H)} - \cot \alpha (H) \frac{1}{\beta (H)} \biggr) \langle \alpha, \beta \rangle \\
& \phantom{aaaaa} + \sum_{\alpha \in \Sigma_r^+}
\frac{m_\alpha(m_\alpha - 2)}{4} \biggl( \text{cosec}^2 \alpha (H) -
\frac{1}{\alpha (H)^2} \biggr) | \alpha |^2 - \sum_{\alpha \in
\Sigma_r^+} \biggl( \frac{m_\alpha |\alpha|}{2} \biggr)^2
\end{align*}

Since the last two terms of the first sum cancel over the summation
for $\alpha \neq 2\beta$, and
$$
\| \rho \|^2 = \langle \rho, \rho \rangle = \sum_{\alpha, \beta \in
\Sigma_r^+} \frac{m_\alpha m_\beta}{4} \langle \alpha, \beta \rangle
= \sum_{\alpha \in \Sigma_r^+} \biggl( \frac{m_\alpha |\alpha|}{2}
\biggr)^2 + \sum_{ \alpha, \beta \in \Sigma_r^+, \; \alpha \neq
\beta} \frac{m_\alpha m_\beta}{4} \langle \alpha, \beta \rangle
$$
we have
\begin{align*}
\Omega_* (H) = -\| \rho \|^2 + 2 \sum_{\alpha \in \Sigma_m^+} \frac{m_\alpha m_{2\alpha}}{4} \biggl( \cot \alpha (H) \cot 2\alpha (H) \biggr) | \alpha |^2 \\
 + \sum_{\alpha \in \Sigma_r^+} \frac{m_\alpha(m_\alpha - 2)}{4} \biggl( \text{cosec}^2 \alpha (H) - \frac{1}{\alpha (H)^2} \biggr) | \alpha |^2
\end{align*}
as claimed. $\phantom{ab} \square$.\\

This expression provides us with an alternate proof of \cite{M},
Lemma 2:
\begin{corollary}  Suppose $X$ is a compact Lie group.  We have
$$
\Omega_* (H) = - \| \rho \|^2
$$
\end{corollary}
\noindent {\bf Proof:} For a compact Lie group, each root (none of
which are multipliable) has multiplicity 2.  Thus it follows that
$F(H) = 0$.
$\; \; \; \square$\\

For compact symmetric spaces that do not have multipliable roots, we
have the following:
\begin{corollary}\label{limit} Suppose $X$ is a compact symmetric space that does not have multipliable
roots, then we have
$$
\lim_{\alpha(H) \rightarrow {\bf 0}} \Omega_* (H) = \sum_{\alpha \in
\Sigma_r^+} \frac{m_\alpha(m_\alpha - 2)}{12} | \alpha |^2
$$
\end{corollary}

\noindent {\bf Proof:} Consider the individual terms of the sum
$$
\sum_{\alpha \in \Sigma_r^+} \frac{m_\alpha(m_\alpha - 2)}{4}
\biggl( \text{cosec}^2 \alpha (H) - \frac{1}{\alpha (H)^2} \biggr) |
\alpha |^2
$$

By l'H\^{o}pital's rule,
$$
\lim_{H \rightarrow H_0} \biggl( \cosec^2 \alpha (H) -
\frac{1}{\alpha (H)^2} \biggr) = \lim_{H \rightarrow H_0} \biggl(
\cosec^2 \alpha (H) - \alpha (H) \cot \alpha (H) \, \cosec^2 \alpha
(H) \biggr)
$$

We re-write the R.H.S. in a neighbourhood of $H_0$ as
\begin{align*}
& \frac{1}{(\alpha (H) - \alpha (H)^3/6 + \mathcal{O}(\alpha (H)^5))^2} \biggl( 1 - \alpha (H) . \frac{1 - \alpha (H)^2/2 + \mathcal{O}(\alpha (H)^4)}{\alpha (H) - \alpha (H)^3/6 + \mathcal{O}(\alpha (H)^4)} \biggr) \\
& \phantom{abcdefg} = \frac{1}{\alpha (H)^2 - \alpha (H)^4/3 + \mathcal{O}(\alpha (H)^6)} \biggl( 1 - \frac{1 - \alpha (H)^2/2 + \mathcal{O}(\alpha (H)^4)}{1 - \alpha (H)^2/6 + \mathcal{O}(\alpha (H)^3)} \biggr)\\
& \phantom{abcdefg}= \frac{1}{\alpha (H)^2(1 - \alpha (H)^2/3 + \mathcal{O}(\alpha (H)^4))} \biggl( \frac{\alpha (H)^2/3 + \mathcal{O}(\alpha (H)^4)}{1 - \alpha (H)^2/6 + \mathcal{O}(\alpha (H)^3)} \biggr)\\
& \phantom{abcdefg}= \frac{1}{3(1 - \alpha (H)^2/3 + \mathcal{O}(\alpha (H)^4))} \biggl( \frac{1+ \mathcal{O}(\alpha (H)^2)}{1 - \alpha (H)^2/6 + \mathcal{O}(\alpha (H)^3)} \biggr) \\
& \phantom{abcdefg}= \frac{1}{3} \biggl( \frac{1 +
\mathcal{O}(\alpha (H)^2)}{1 - \mathcal{O}(\alpha (H)^2)} \biggr)
\longrightarrow \frac{1}{3} \; \; \text{as} \; \; H \longrightarrow
H_0.
\end{align*}
The proposition follows. $\phantom{abc} \square$\\

From the classification of symmetric spaces (see \cite{HE1}, Ch X),
the only simple, simply connected compact symmetric spaces of rank
one are the $n$-dimensional spheres, $S^n$.  For $S^n$ there is one
root $\alpha$ which has multiplicity $n-1$. Normalising $\alpha =
1$, we have the following:
\begin{corollary}\label{Slimit}  Suppose $X = S^n$.  We have
$$
\Omega_* (H) = -\biggl( \frac{n-1}{2} \biggr)^2 + \biggl(
\frac{(n-1)(n-3)}{4} \biggr) \biggl( \cosec^2 H - \frac{1}{H^2}
\biggr)
$$

Furthermore,
$$
\lim_{H \rightarrow 0} \Omega_* (H) = \tfrac 16 n (1-n)
$$
\end{corollary}

However, we have not been able to calculate the heat kernel with
these potentials on $\p \cong \R^n$.  Despite this, we have the
following theorem guaranteeing the existence of a fundamental
solution whose wrap is the heat kernel on $X$:
\begin{theorem}\label{HorThm}
There exists a fundamental solution of the semigroup $\exp t(L_\p +
\Omega_*)$ in the fundamental domain of $X$ whose wrap is the heat
kernel on $X$.
\end{theorem}

{\bf Proof:} By H\"ormander's theorem,
$$
\frac{\partial}{\partial t} - (L_\p + \Omega_*)
$$
is hypoelliptic on the set $\{ H \, : \, |\alpha(H)|< \pi -
\epsilon, \; \forall \alpha \in \Sigma^+_r \; \text{and} \; \epsilon
> 0 \}$, that is, the fundamental domain of $X$.  Therefore, $\tfrac{\partial}{\partial t} - (L_\p +
\Omega_*)$ has a fundamental solution in the fundamental domain of
$X$.  This fundamental solution may be wrapped to $X$ (since it only
takes values in the fundamental domain of $X$), which by Theorem
\ref{peturbHK2} is the heat kernel on $X$. $\phantom{abc} \square$

\subsection{``Radial'' wrapping of the Laplacian}

We now consider an alternate approach that may allow us to verify
our results by computing the wrap of the Laplacian for $U/K$ in
spherical polar co-ordinates.  This has the effect of breaking up
the potential $j^{-1} L_\p j$ into two separate terms on $\p$ and on
$U/K$ where there is some literature on heat kernels with these
potentials. However, the known results are not sufficient to verify our conjectures.\\

We will now prove a general result regarding the wrapping the
Laplacian on compact symmetric spaces.  As we stated in the first
half of the thesis, the Laplacian is of fundamental importance as it
is the infinitesimal generator of the heat semigroup and Brownian
motion.  Since we are only considering the wrap of $K$-invariant
functions on $\p$ to $K$-invariant functions on $U/K$, it will
suffice to consider the radial parts of these functions on $\p$ and
$U/K$.\\

Using the facts that the wrap is linear and the Laplacian may be
decomposed into radial and angular parts, $\Delta (L_{\p})$ and
$L^A_{\p}$, respectively, we have
\begin{align*}
\langle \Phi(L_{\p} u), f \rangle & = \langle  (\Delta (L_{\p}) + L^A_{\p}) u , j \tilde{f}\rangle\\
& = \langle  \Delta (L_{\p})u, j \tilde{f}) + \langle L^A_{\p}u, j
\tilde{f}\rangle
\end{align*}
so that the wrap of radial and angular parts may be considered separately.\\

Firstly, conformality of the exponential map ensures that
transversal manifolds are the same on the symmetric spaces and their
tangent spaces under $\Exp$.  Conformality of $\Exp$ also ensures
that the Laplacians on transversal manifolds under the exponential
map are equivalent (see \cite{HE2}, Ch. 2, Thm. 3.15).\\

In general, the radial part of the Laplacian may be given by the
following:
\begin{theorem}\label{631} (c.f. \cite{HE2}, Ch. 2, Thm. 3.7) Let $V$ be a Riemannian manifold, $H$ a closed unimodular subgroup of the isometry group $I(V)$.  Assume that a submanifold $W \subset V$ satisfies the following transversality condition:  For each $w \in W$,
$$
(H \cdot w) \cap W = \{w\}, \phantom{abcd} V_w = (H \cdot w)_w
\oplus W_w
$$
(orthogonal direct sum).  Then the radial part of $L_V$ is given by
$$
\Delta (L_V) = \delta^{-1/2} L_W \circ \delta^{1/2} - \delta^{-1/2}
L_W (\delta^{1/2} )
$$
where $\circ$ denotes the composition of operators and $\delta$ the
density function.
\end{theorem}

{\bf Example:}  For the case of $V = \R^n$ under the operation of $H
= O(n)$ (c.f. \cite{HE2}, Ch. 2, pp 266), the submanifold $W= \R^+ -
\{0\}$ satisfies the transversality condition, and $\delta =
c|x|^{n-1}$.  On $\R^2$ this simplifies to
\begin{equation}
\Delta (L_{\R^2}) = r^{-1/2} \frac{\partial^2}{\partial r^2} \circ
r^{1/2} + \tfrac 14 r^{-2}
\end{equation}

We will keep our notation consistent with \cite{HE2}, page 273, by
letting
$$
\delta_0 (H) = \prod_{\alpha \in \Sigma_r^+} \alpha (H)^{m_\alpha}
$$
be the density function for $K$ acting on $\p \cong \R^n$, and
$$
\delta (\Exp H) = \prod_{\alpha \in \Sigma_r^+} (\sin \alpha
(H))^{m_\alpha}
$$
be the density function for $K$ acting on $U/K$, with the Jacobian
of $\Exp$ being $J = \delta / \delta_0$.  For a general compact
symmetric space, we now have the following theorem showing how the
radial parts of the Laplacian on $\p$ wrap to the radial parts of
$U/K$:
\begin{theorem}\label{wraplapcss}  With the above notation, the radial part of $L_{\p}$, given by
\begin{equation}\label{lapp}
\Delta (L_{\p}) = \delta_0^{-1/2} L_\af \circ \delta_0^{1/2} -
\delta_0^{-1/2} L_\af (\delta_0^{1/2} ),
\end{equation}
wraps to the radial part of $L_{U/K}$, given by
\begin{equation}\label{lapuk}
\Delta (L_{U/K}) = \delta^{-1/2} L_A \circ \delta^{1/2} -
\delta^{-1/2} L_A (\delta^{1/2} ).
\end{equation}

That is,
$$
\Phi (\Delta (L_{\p})) = (\Delta (L_{U/K})) \Phi
$$
\end{theorem}

\noindent {\bf Proof:}  By using the fact that the Laplacian is a
symmetric operator, Theorem \ref{631}, and \cite{HE2} Ch. II, $\S
3$, we have,
\begin{align*}
\bigl\langle \Phi (\Delta (L_{\p}) u ) , f \bigr\rangle & =
\bigl\langle (\Delta (L_{\p})) u, j \tilde{f} \bigr\rangle \tag{by
definition of the
wrap}\\
& = \bigl\langle u, (\Delta (L_{\p})) (j \tilde{f}) \bigr\rangle
\tag{Since $\Delta (L_{\p})$ is a symmetric operator}
\end{align*}

By Theorem \ref{631} this is
$$
\bigl\langle \Phi (\Delta (L_{\p}) u ) , f \bigr\rangle =
\bigl\langle u, j \, j^{-1} \bigl( \delta_0^{-1/2} L_\af \circ
\delta_0^{1/2} - \delta_0^{-1/2} L_\af (\delta_0^{1/2} ) \bigr) (j
\tilde{f}) \bigr\rangle
$$

From the proof of \cite{HE2} Ch. II, Thm. 3.15, this equal to
\begin{align*}
\bigl\langle \Phi (\Delta (L_{\p}) u ) , f \bigr\rangle & = \bigl\langle u, j \, \bigl( \delta^{-1/2} L_A \circ \delta^{1/2} - \delta^{-1/2} L_A (\delta^{1/2} ) \bigr) \tilde{f} \bigr\rangle\\
& = \bigl\langle u, j \, \widetilde{(\Delta (L_{U/K})f)} \bigr\rangle\\
& = \bigl\langle \Phi u, (\Delta (L_{U/K})f \bigr\rangle\\
& = \bigl\langle (\Delta (L_{U/K})) \Phi u  , f \bigr\rangle
\end{align*}
as required.  $\phantom{abc} \square$\\

We will now show explicitly the details of the above calculation for
the 2-sphere.  Consider the decomposition of the Laplacian into its
radial and angular parts on both the sphere and the plane. Using
spherical polar co-ordinates we have:
$$
x = \cos \psi \sin \theta, \; \; y = \sin \psi \sin \theta, \; \; z
= \cos \theta
$$
for $\psi \in [0, 2\pi)$ and $\theta \in [0, \pi)$.\\

The Laplacian on $S^2$ is given by
\begin{equation}
L_{S^2} = \frac{\partial^2}{\partial \theta^2} + \cot \theta
\frac{\partial}{\partial \theta} + \sin^{-2} \theta
\frac{\partial^2}{\partial \psi^2}
\end{equation}

Now, the Laplacian on $\R^2$, given in polar co-ordinates is
\begin{equation}
L_{\R^2} = \frac{\partial^2}{\partial r^2} + \frac{1}{r}
\frac{\partial}{\partial r} + \frac{1}{r^2}
\frac{\partial^2}{\partial \phi^2}
\end{equation}
and note the polar decomposition into radial and angular parts,
where the angular part is given by the last term, and may be
considered as the Laplacian on the circle of radius $r$, ie,
\begin{equation}
L = \frac{1}{r^2} \frac{\partial^2}{\partial \phi^2}
\end{equation}

Now consider the wrapping map:  Let $v$ be a $K$-invariant
distribution of compact support or Schwartz function on $\R^2$.  The
wrapping map, $\Phi$ is given by
\begin{equation}
\langle \Phi(v), f \rangle = \langle v,j \tilde{f} \rangle
\end{equation}

The $j$-function in this case is $j(X) = \Bigl( \frac{\sin X}{X}
\Bigr)^{1/2}$.  Essentially, we are using the exponential map which
will map the variables to one another, ie,
\begin{align*}
\phi & \longrightarrow \psi\\
r & \longrightarrow \theta
\end{align*}
which we may consider independently since the angular and radial
parts commute.  The wrap contains information about the curvature.
The map $\phi \rightarrow \psi$ has no curvature (ie, $j(X) = 1$),
but $r \longrightarrow \theta$ does.\\

We will now give the calculations for Theorem \ref{wraplapcss} in
the case of the 2-sphere.  The radial part of the Laplacian on
$\R^2$ is
\begin{equation}
\Delta (L_{\R^2}) = \frac{\partial^2}{\partial r^2} + \frac{1}{r}
\frac{\partial}{\partial r}
\end{equation}
which for may be expressed as
\begin{equation}
\Delta (L_{\R^2}) = r^{-1/2} \frac{\partial^2}{\partial r^2} \circ
r^{1/2} + \tfrac 14 r^{-2}
\end{equation}

The rays of the radial part are being wrapped to the meridians of
longitude on the sphere, weighted by the $j$-function.  The density
function is given by
\begin{equation}
\delta (\theta = \Exp r)^{1/2} = \prod_{\alpha \in \Sigma_r^+} \sin
\alpha (\theta = \Exp r) = j^2 (r) \prod_{\alpha \in \Sigma_r^+}
\alpha(r)
\end{equation}

For the $2$-sphere there is only one root with multiplicity 1, and
so $\delta (\theta = \Exp r)^{1/2} = \sin (\theta)$.  Now,
\begin{align*}
\Delta (L_{\R^2})f & = r^{-1/2} \frac{\partial^2}{\partial r^2} \circ r^{1/2} f + \tfrac 14 r^{-2}f\\
& = r^{-1/2} \biggl( \frac{\partial^2}{\partial r^2} \circ r^{1/2}f \biggr) + \tfrac 14 r^{-2}f\\
& = r^{-1/2} \biggl( \frac{\partial^2}{\partial r^2} (r^{1/2}f) \biggr) + \tfrac 14 r^{-2}f\\
& = r^{-1/2} \biggl( \tfrac 14 r^{-3/2}f + r^{-1/2}\frac{\partial f}{\partial r} + r^{1/2} \frac{\partial^2 f}{\partial r^2} \biggr) + \tfrac 14 r^{-2}f\\
& = r^{-1}\frac{\partial f}{\partial r} + \frac{\partial^2
f}{\partial r^2}
\end{align*}
so that
$$
\Delta (L_{\R^2}) = \frac{1}{r} \frac{\partial}{\partial r} +
\frac{\partial^2}{\partial r^2}
$$
as required.\\

For the two-sphere we have $\delta (\theta = \exp r)^{1/2} = \sin
(\theta)$, since $r$ becomes $\sin \theta$ under the wrapping map on
$L_{S^2}$. The angular part of the Laplacians on the tranversal
manifolds remains essentially the same, changing only from $L_W =
\frac{\partial^2}{\partial r^2}$ to $L_W =
\frac{\partial^2}{\partial \theta^2}$.  By Theorem \ref{631} the
radial part of $L_{S^2}$ is
$$
\Delta (L_{S^2}) = (\sin \theta)^{-1/2} L_W \circ (\sin
\theta)^{1/2} - (\sin \theta)^{-1/2} L_W ((\sin \theta)^{1/2} ).
$$

Firstly, let
$$
S = \frac{\partial^2}{\partial \theta^2} (\sin \theta)^{1/2} =
-\tfrac 14 \cos^2 \theta (\sin \theta)^{-3/2} - \tfrac 12 (\sin
\theta)^{1/2}
$$
\begin{align*}
\Phi\bigl( \Delta (L_{\R^2})f \bigr) & = (\sin \theta)^{-1/2} \frac{\partial^2}{\partial \theta^2} \circ (\sin \theta)^{1/2} f - (\sin \theta)^{-1/2}Sf\\
& = (\sin \theta)^{-1/2} \biggl( Sf + 2.\tfrac 12 \cos \theta (\sin \theta)^{-1/2} \frac{\partial}{\partial \theta}f + (\sin \theta)^{1/2} \frac{\partial^2}{\partial \theta^2}f  \biggr) - (\sin \theta)^{-1/2}S f\\
& = \cos \theta (\sin \theta)^{-1} \frac{\partial}{\partial \theta}f + \frac{\partial^2}{\partial \theta^2}f\\
& = \cot \theta \frac{\partial f}{\partial \theta} +
\frac{\partial^2 f}{\partial \theta^2}
\end{align*}
so that
$$
\Phi\bigl( \Delta (L_{\R^2})f \bigr) = \frac{\partial^2}{\partial
\theta^2} + \cot\theta \frac{\partial}{\partial \theta}
$$
as required.\\

The potential terms on $\p$ and $U/K$, $\delta_0^{-1/2} L_\af
(\delta_0^{1/2} )$ and $\delta^{-1/2} L_A (\delta^{1/2} )$,
respectively are somewhat manageable.\\

For the case of the two-sphere, the potential for the radial
Laplacian on $S^2$ is
$$
\delta^{-1/2} L_A (\delta^{1/2} ) = (\sin H)^{-1/2}
\frac{\partial}{\partial H} (\sin H)^{1/2} = \frac{1}{4} -\frac{1}{4
\sin^2 (\theta)}
$$
and the potential for the radial Laplacian on $\R^2$ is
$$
\delta^{-1/2} L_{\af} (\delta^{1/2} ) = H^{-1/2}
\frac{\partial}{\partial H} H^{1/2} = \frac 14 \frac{1}{H^2}
$$

Thus, it remains to find the solution to the radial heat equation in
potential form, that is, the solution to
\begin{equation}\label{PotSph}
\frac{\partial u}{\partial t} = \Bigl[ \frac{1}{2}
\frac{\partial^2}{\partial \theta^2} -\frac 18 + \frac 18 \cosec^2
\theta \Bigr] (u)
\end{equation}
on $S^2$, and
\begin{equation}\label{PotSph2}
\frac{\partial u}{\partial t} = \Bigl[ \frac{1}{2}
\frac{\partial^2}{\partial \theta^2} -\frac 18 \frac{1}{H^2} \Bigr]
(u)
\end{equation}
on $\R^2$.\\

The solution to the problem of this heat equation with the potential
in (\ref{PotSph}) above has been investigated by Peak and Inomata
(\cite{PI}); see also Anderson and Anderson (\cite{AA}).  The
problem of Brownian motion on the 2-sphere is also known as the
``rigid rotor'' problem. The approach of \cite{PI} is to apply a
constraint to a Brownian motion in $\R^3$, leading to the solution
to a certain path integral.  This has been generalised by Anderson
(\cite{A}) to the general $n$-sphere by intertwining the Laplacian
on the 2-sphere with a certain fractional differential operator. The
heat kernel on $S^n$ is then found by applying this operator to the
heat kernel on the $2$-sphere.  However, generalising this approach
to rank 2 spaces and beyond seems quite difficult (\cite{AC}).\\

To find the heat kernel via wrapping, we would need to find the heat
kernel with potential $1/H^2$ on $\p$ and wrap it to $U/K$.
Existence results exist for the heat kernel with this potential (see
\cite{BG}), but we have been unable to find an explicit form for it.

\section{The Non-compact Case}

In this section we consider the wrapping map, Brownian motion and
heat kernels on certain non-compact symmetric spaces, namely complex
Lie groups $G_\C$, the symmetric spaces $G_\C /K$, and the symmetric
spaces of split rank.\\

We extend the wrapping map to complex Lie groups, and then use this
to wrap Brownian motion and the heat kernel on these spaces. Our
results also hold for the symmetric spaces $G_\C /K$, although for
this case we have not been able to prove a global wrapping formula.
Finally, we conclude by considering the symmetric spaces of split
rank, and deduce a similar result to that of above, in that
``bending" the heat kernel from the tangent space does not yield the
heat kernel on $G/K$ since the convolution structure is not
preserved.\\

Firstly, we state the quantity $j^{-1} L_\p j$ for the non-compact
case.  This calculation is essentially the same as that for the
quantity $(j \circ \log) L_{G/K} (j \circ \log)^{-1}$ given in
\cite{HS}, section 2.

\begin{prop}\label{prop721}  With the above notation,
$$
(j^{-1} L_\p j) (H) = \| \rho \|^2 + F(H)
$$
where
\begin{align*}
F(H) & = \phantom{aaa} \sum_{\alpha \in \Sigma_r^+} \frac{m_\alpha(m_\alpha - 2)}{4} \biggl( \csch^2 \alpha (H) - \frac{1}{\alpha (H)^2} \biggr) | \alpha |^2\\
& + 2 \sum_{ \alpha \in \Sigma_m^+} \frac{m_\alpha m_{2\alpha}}{4} \biggl( \csch^2 \alpha (H) - \frac{1}{\alpha(H)^2}  \biggr) |\alpha|^2\\
\end{align*}
\end{prop}

When $X$ is a complex Lie group, or a symmetric space $G/K$, $G$
complex, each root has a multiplicity of two, none of which are
multipliable.  It is easily seen that $F(H) = 0$, and we
consequently have:
\begin{corollary}\label{722}  Suppose $X$ is a complex Lie group, or a symmetric space $G/K$, $G$ complex, then
$$
(j^{-1} L_\p j) (H) = \| \rho \|^2
$$
\end{corollary}

\subsection{Complex Lie Groups}

We briefly consolidate some notations and definitions of complex Lie
groups. Let $G_\C$ be a complex, connected, semisimple Lie group
with Lie algebra $\g_\C$.  Let $\g_\C = \kg + \p$ be a Cartan
decomposition of $\g_\C$, and $K$ the compact group corresponding to
$\kg$. We denote by $\g^\R$ as $\g_\C$ realised as a real Lie
algebra, with the Cartan decomposition given by $\g_\C = \kg +
i\kg$, and $\kg$ may be identified as the compact real form.\\

We let $\af \subset \p$ be a maximal abelian subalgebra of $\p$, and
$\Sigma_r = \Sigma_r (\g_\C,\af)$ the root system of $\g_\C$ with
respect to $\af$.  The real dimension of each root space is two,
that is $m_\alpha = 2$ for all $\alpha \in \Sigma_r$.  Elements of
$G_\C$ may be written as $k_1 a k_2$, where $k_1, k_2 \in K$ and $a
\in \exp (\bar{\af}^+)$.  Since the roots are purely real on $\af$,
$A^+ = \exp (\bar{\af}^+)$ is of the form $(\R^+)^l$. Furthermore,
from \cite{KNP} Ch. II, Thm. 2.15, any two Cartan subalgebras are
conjugate under the adjoint action of
$G_\C$.\\

A real-valued square root of the Jacobian of the exponential map is
given by
$$
j(H) = \prod_{\alpha \in \Sigma_r^+} \biggl( \frac{\sinh \alpha
(H)/2}{\alpha (H)/2} \biggr), \phantom{abcde} H \in \af^+
$$

Since each root of a complex Lie group has multiplicity 2 we have
(see \cite{HE3}, pp 487) $j^{-1} L j = \| \rho \|^2 = \frac{\dim
G}{12}$. (Compare this to the case of a compact Lie group, where we
have $ -j^{-1} L j = \| \rho \|^2 = \frac{\dim G}{24}$.  The
elementary spherical functions on a complex Lie groups are given by
(see \cite{HE2}, Ch. IV, Thm. 5.7):
$$
\varphi_\lambda (a) = c(\lambda) \frac{\sum_{\omega \in W} \sgn
\omega \, e^{i \omega \lambda (\log a)} }{\sum_{\omega \in W} \sgn
\omega \, e^{\omega \rho (\log a)} }, \phantom{abcde} a \in A
$$
where the function $c$ is given by
$$
c(\lambda) = \prod_{\alpha \in \Sigma_r^+} \langle \alpha , \rho
\rangle / \prod_{\alpha \in \Sigma_r^+} \langle \alpha ,\lambda
\rangle
$$

We also note another formula of Harish-Chandra (see \cite{HE2}, Ch.
II, Thm. 5.35):  Let $U = G/K$ and let $du$ be normalised Haar
measure on $U$, and let $\pi$ be the product of positive roots.
Then, if $H, \, H' \in \tg_\C$,
\begin{equation}\label{733}
\pi(H) \pi(H') \int_U e^{\langle u H, H' \rangle} du =
\tfrac{1}{|W|} \partial(\pi)(\pi) \sum_{\omega \in W} \sgn \omega \,
e^{\langle \omega H, H' \rangle}
\end{equation}

Setting $H' = H_\rho$, we obtain (\cite{HE2}, Ch. II, Cor. 5.36)
$$
|W| = \partial (\pi)(\pi) / \pi (H_\rho)
$$

The heat kernel of a complex group is given by the following:
\begin{prop}\label{HKComplexLG} (\cite{GAN}) Suppose $G_\C$ is a complex Lie group.  The heat kernel on $G_\C$, $q_t (a)$, given on $A$ is
\begin{equation}
q_t (a) = (2\pi t)^{-n/2} \exp (-t \| \rho \|^2 /2) \frac{1}{j(\log
a)} \exp (-|\log a|^2 /2t), \phantom{abcde} a \in A
\end{equation}
\end{prop}

We now prove the wrapping theorem for complex groups.  The wrapping
map $\Phi$ is given by
$$
\langle \Phi(\nu), f \rangle = \langle \nu , j \tilde{f} \rangle
$$
and is well defined if $\nu$ is an integrable $K$-bi-invariant function.\\

We face a difficulty with the orbital convolution theory, since the
co-adjoint orbits are non-compact, and are not in one-to-one
correspondence with the adjoint orbits, as in the compact case.
However, the wrapping formula for complex groups does hold for
spherical measures, since the Fourier theory is essentially reduced
to the abelian case.
\begin{theorem} Suppose $\mu, \nu$ are two $K$-bi-invariant Schwartz functions on $\g_\C$, then
$$
\Phi (\mu * \nu) = \Phi(\mu) * \Phi(\nu)
$$
\end{theorem}
{\bf Proof:} The spherical Fourier transform of $\Phi(u)$ at a
representation $\pi$, indexed by highest weight $\lambda$, is given
by
\begin{align*}
\langle \Phi (u), \varphi_\lambda \rangle & = \langle u, j \tilde{\varphi}_\lambda \rangle\\
& = \Bigl\langle u, c(\lambda) \frac{\sum_{\omega \in W} \sgn \omega \, e^{i \omega \lambda (\cdot)}}{\prod_{\alpha > 0} \alpha (H) } \Bigr\rangle\\
& = \frac{\pi(\rho)}{\pi(i\lambda)} \Bigl\langle u,
\frac{\sum_{\omega \in W} \sgn \omega \, e^{i \omega \lambda
(\cdot)}}{\prod_{\alpha > 0} \alpha (H) } \Bigr\rangle
\end{align*}

By Harish-Chandra's formula, this is
\begin{align*}
\langle \Phi (u), \varphi_\lambda \rangle & = \Bigr\langle u, \int_K e^{i k \lambda (\cdot) } dk \Bigr\rangle\\
& = \langle u, e^{i\lambda (\cdot)} \rangle \tag{Since $u$ is $K$-bi-invariant}\\
& = \widehat{u} (\lambda).
\end{align*}
Using abelian Fourier analysis, we obtain
\begin{align*}
(\Phi(\mu * \nu))^\wedge & = (\mu * \nu)^\wedge{}\\
& = \mu^\wedge \cdot \nu^\wedge\\
& = \Phi^\wedge (\mu) \cdot \Phi^\wedge (\nu)\\
& = (\Phi (\mu) * \Phi (\nu))^\wedge. \phantom{abcde} \square
\end{align*}



\noindent {\bf Remark:} For Lie groups, Rouvi\`ere's conjecture that
$e=1$ for Lie groups is equivalent to the Kashiwara-Vergne
conjecture (\cite{KV}).  This was recently proven in \cite{AST}. For
the symmetric spaces $G_\C / K$, $e=1$ was proven in a neighbourhood
of $0 \in \s$ by Torossian in \cite{T}.

\subsection{The wrap of Brownian motion and heat kernels on complex
Lie groups}

We prove the following general result for computing the wrap of a
$K$-bi-invariant Schwartz function:
\begin{prop}\label{prop742}  Let $G/K$ be a semisimple Riemannian symmetric space of the non-compact type, with tangent space $\p$.  Let $\varphi \in S(\p)$ be $K$-bi-invariant.  Then $\Phi (j \varphi)$ is a $C^\infty$ $K$-bi-invariant function, given on $A$,
by
$$
\Phi (j \varphi) (\exp H) = \varphi (H), \phantom{abcde} H \in \af.
$$
\end{prop}

\noindent {\bf Proof:} Firstly, we let
$$
f^K(g) = \int \int_{K \times K} f(k_1 g k_2) dk_1 dk_2
$$

Since $\Exp : \p \rightarrow G/K$ is a global diffeomorphism, it
follows that for any $f \in C^\infty (G/K)$ that
\begin{align}
\langle \Phi (j \varphi) , f \rangle & = \langle j \varphi , j \tilde{f} \rangle = \int_\p j^2 (X) \varphi (X) \tilde{f} (X) dX\label{742}\\
& = \frac{1}{|W|} \int_{\af} \prod_{\alpha \in \Sigma_r^+} (\sinh \alpha (H)/2) \varphi (H) \tilde{f}^K (H) dH\label{743}\\
& = \frac{1}{|W|} \int_{A} \prod_{\alpha \in \Sigma_r^+} (\sinh \alpha (\log a)/2) \varphi (\log a) \tilde{f}^K (a) da\label{744}\\
& = \int_{G/K} \varphi (\log g) f(g) dg \label{745}
\end{align}

The result follows from \ref{5214} $\phantom{abc} \square$\\

Using Proposition \ref{prop742} for complex Lie groups, we
immediately recover the results of Proposition \ref{HKComplexLG}:
\begin{corollary}\label{cor743}  Let $p_t (H) = (2\pi t)^{-n/2} e^{-|H|^2 / 2t}$ be the heat kernel on $\g$ (given on $\af$).  Then
$$
\Phi (p_t) (\exp H) = (2\pi t)^{-n/2} \frac{1}{j(H)} \exp (-|H|^2
/2t), \phantom{abcde} H \in \af
$$
\end{corollary}
which is the shifted heat kernel, $q_t^\rho (g)$.\\

\noindent {\bf Remark:}  At this point we identify two critical
steps in the proof of Proposition \ref{prop742}, namely (\ref{743})
to (\ref{744}), and (\ref{744}) to (\ref{745}).  In \cite{DW2}, the
step (\ref{743}) to (\ref{744}) is achieved for a compact Lie group
through a function $\Psi$, which is a sum over the integer lattice
of the function $\varphi$.  Thus (\ref{743}) may be given on a
fundamental domain, and hence (\ref{744}) obtained.  The step in
going from (\ref{744}) to (\ref{745}) relies on the fact that all
Cartan subalgebras in $\g$ are conjugate to each other under the
adjoint group.\\

We now consider the wrap of the Laplacian:
\begin{prop}\label{prop744}  Let $G$ be a complex connected Lie group with Lie algebra $\g$.  Then for any $u \in S(\g)$
$$
\Phi \bigl( L_\g (u) \bigr) = (L_G + \| \rho \|^2) \bigl( \Phi u
\bigr)
$$
where $\Phi$ is the wrapping map, $L_\g$ is the Laplacian on $\g$
(regarded as a Euclidean vector space), and $\rho$ the half sum of
positive roots.
\end{prop}

We now wrap Brownian motion as previously, by defining shifted
Brownian motion $\xi_t$ on $G$, and lifting $f$ in the It\^o
equation by $f \mapsto j.f \circ \exp = h$. That is,
\begin{df}\label{prop745}  Let $\zeta_t$ be a Brownian motion on $\g \cong \R^n$.  The wrap of $\zeta_t$ is to a process $\xi_t$ on $G$ is given by the mapping $f \mapsto j.f \circ \exp = h$, where $\xi_t$ is shifted Brownian motion.  We write this as
$$
\Phi (\zeta_t) = \xi_t
$$
\end{df}

From this we have the analogue of our result from \cite{M2},
regarding the wrap of Brownian motion:
\begin{theorem}\label{thm746}
Suppose $\xi_t$ is the wrap of the Brownian motion on $\g$,
$\zeta_t$.  Then the law of $\xi_t$ may be found by wrapping the law
of Brownian motion on its Lie algebra.  That is,
$$
\E_X(j. f \circ \exp(\zeta_t)) = \E_{\exp X} (f(\xi_t))
$$
which in law is given by
$$\Phi (p_t)(\exp H) = q_t^\rho (g)$$
where $p_t(x)$ is the heat kernel on $\g = \R^n$, and $q_t^\rho (g)$
is the heat kernel corresponding to the shifted Laplacian on $G$
\end{theorem}

As in \cite{M2}, we apply the Feynman-Ka\v{c} formula to obtain:
\begin{prop}\label{prop746}
The expection of $(\xi_t)_{t \geq 0}$ under $\tilde{\bbP}$, denoted
by $\tilde{\E}$, is
$$\tilde{\E}(\xi_t) = (2\pi t)^{-n/2} e^{-\| \rho \|^2 t/2} \frac{1}{j(H)} \exp (-|H|^2 /2t), \phantom{abc} t \in \R^+, \; H \in \af.$$
which is the standard heat kernel on $G$.
\end{prop}
{\bf Proof:}  Taking expectations under $\tilde{\bbP}$ yields
\begin{align*}
\tilde{\E}(\xi_t) & = \E(e^{-\| \rho \|^2 t/2} \xi_t) = e^{-\| \rho \|^2 t/2} \E(\xi_t) = e^{-\| \rho \|^2 t/2}q_t^\rho (g)\\
& = q_t (g) = (2\pi t)^{-n/2} e^{-\| \rho \|^2 t/2} \frac{1}{j(H)}
\exp (-|H|^2 /2t), \phantom{abc} t \in \R^+, \; H \in \af.
\end{align*}
as required. $\phantom{abcde} \square$

\subsection{Wrapping for certain non-compact symmetric spaces}

There are many examples of non-compact symmetric spaces.  In this
section we conjecture that our results on wrapping Brownian motion
and heat kernel may be extended to other symmetric spaces  -  more
specifically, our results on wrapping Brownian motion and heat
kernels should hold for every symmetric space for which the wrapping
theorem holds.  The wrapping theorem is conjectured to
hold for all Lie groups in \cite{DW2}.\\

Firstly, we will consider the symmetric spaces $G/K$, $G$ complex.
The wrapping theorem was proved to hold locally for these spaces in
\cite{T}. We have not been able to prove a global wrapping theorem,
but analogous results of section 7.5 hold.\\

Recall from Corollary \ref{722} that $(j^{-1} L_\p j) = \| \rho
\|^2$.  We apply Proposition \ref{prop742} to find the shifted heat
kernel, and wrap Brownian motion and the heat kernel as we did in
the complex case, yielding the following analogue of Corollary
\ref{cor743}:
\begin{prop}\label{prop756}
The expectation of $(\xi_t)_{t \geq 0}$ under $\tilde{\bbP}$,
denoted by $\tilde{\E}$, is
$$
\tilde{\E}(\xi_t) = (2\pi t)^{-n/2} e^{-\| \rho \|^2 t/2}
\frac{1}{j(H)} \exp (-|H|^2 /2t), \phantom{abc} t \in \R^+, \; H \in
\af.
$$
which is the standard heat kernel on $G$.
\end{prop}

{\bf Remark:}  Proposition \ref{prop756} holds with $\| \rho \|^2 =
\tfrac{n}{12}$, which agrees with the expression for the heat kernel
on the symmetric spaces $G/K$,
$G$ complex, given by Arede in \cite{ARE}, $\S 2.2$.\\

We now consider the symmetric spaces of {\it split-rank} type:
\begin{df}  A symmetric space $G/K$ is said to be of {\bf split-rank} if {\rm rank} $G$ = {\rm rank} $G/K$ + {\rm rank} $K$.
\end{df}

\begin{prop}  (\cite{HE1}, Ch. IX, Thm. 6.1)  The following are equivalent:
\begin{list}{}{}
\item (i)  $G/K$ has split rank,
\item (ii)  Each restricted root has even multiplicity,
\item (iii)  All Cartan subalgebras of $\g$ are conjugate under the adjoint action.
\end{list}
\end{prop}

It follows from (ii) and (iii) that the roots of $G/K$ with respect
to $\af_\C$ are real on $\af$.  Thus, the maximal abelian subgroup
corresponding to the Cartan subalgebra is of the form $(e^{\alpha_1
(H)}, \dots , e^{\alpha_l (H)})$, where the $\alpha_i$'s are all
real, so the subgroup is of the form $(\R^+)^n$.\\

By \cite{HE1}, Ch. IX, \S 4, it is known that the non-compact,
simple, simply connected symmetric spaces $X = G/K$ of split-rank
consist of the following:
\begin{list}{}{}
\item (i)  $\R^{2n} \rtimes \R^+ \cong SO_0 (2n+1,1) / SO(2n+1)$, the odd dimensional hyperbolic spaces,
\item (ii)  $G_\C / K$, where $K$ is a maximal compact subgroup of $G_\C$,
\item (iii)  $SU^*(2n) / Sp (n)$,
\item (iv)  $E_{6 (-26)}/ F_4$.
\end{list}

We apply proposition \ref{prop742} for these spaces, and the heat
kernel $p_t$ may be wrapped from $\p$ to give
\begin{equation}\label{GANC}
\Phi(p_t)(\Exp H) = (2\pi t)^{-n/2} e^{-\| \rho \|^2 t/2}
\frac{1}{j(H)} \exp (-|H|^2 /2t), \phantom{abc} t \in \R^+, \; H \in
\af.
\end{equation}

However, the $e$-function is not identically 1 for all the spaces of
split-rank type.  We now make a similar observation to to that for
compact symmetric spaces: (\ref{GANC}) is not the heat kernel for
the spaces of split-rank type when the $e$-function is not
identically 1. (\ref{GANC}) is the Gaussian approximation to the
heat kernel, and formally state:
\begin{theorem} The Gaussian approximation (\ref{GANC}) is not
exactly equal to the heat kernel (modulo the phase factor of $e^{\|
\rho \|^2 t}$) for the non-compact symmetric spaces of split-rank
type that have a non-trivial $e$-functions.
\end{theorem}

For non-compact symmetric spaces that do not have multipliable
roots, we have the following:
\begin{prop} Suppose $X$ is a non-compact symmetric space that does not have multipliable
roots, then $\Omega_*$ is a bounded $C^\infty$ function on $\p$.
Furthermore,
$$
\lim_{\alpha(H) \rightarrow {\bf 0}} \Omega_* (H) = -\sum_{\alpha
\in \Sigma_r^+} \frac{m_\alpha(m_\alpha - 2)}{12} | \alpha |^2
$$
\end{prop}

\noindent {\bf Proof:} It is clear that $\Omega_*$ is a bounded
$C^\infty$ function on $\p$, except on the hyperplane in $\p$ where
$\alpha (H) = 0$, which we will denote $H_0$. To show that
$\Omega_*$ is bounded on $H_0$, we consider the individual terms of
the sum
$$
\sum_{\alpha \in \Sigma_r^+} \frac{m_\alpha(m_\alpha - 2)}{4}
\biggl( \csch^2 \alpha (H) - \frac{1}{\alpha (H)^2} \biggr) | \alpha
|^2
$$

The proof is now analogous to Corollary \ref{limit}, replacing cosec by csch, and noting the change in sign. $\phantom{abc} \square$\\

From the classification of symmetric spaces (see \cite{HE1}, Ch X),
the only simple, simply connected non-compact symmetric spaces of
rank one with no multipliable roots are the $n$-dimensional
hyperbolic spaces, $H^n \cong SO_0 (n,1) / SO(n)$, $n \geq 2$. Since
$H^n$, is a symmetric space of rank 1, $\dim \af = 1$ and so $\af
\cong \R$. Therefore there is only one root $\alpha \in \af^*$ which
has multiplicity $n-1$.  We normalise $\alpha$ such that $\alpha =
1$, thus $\rho = (n-1)/2$, giving us the non-compact analogue of
proposition \ref{Slimit}:
\begin{corollary}  Suppose $X = H^n$.  We have
$$
\Omega_* (H) = \biggl( \frac{n-1}{2} \biggr)^2 + \biggl(
\frac{(n-1)(n-3)}{4} \biggr) \biggl( \csch^2 H - \frac{1}{H^2}
\biggr)
$$

Furthermore,
$$
\lim_{H \rightarrow 0} \Omega_* (H) = \tfrac 13 (2n^n +5n +3)
$$
\end{corollary}

The $e$-function for $H^n \cong SO_0 (n,1) / SO(n)$, $n \geq 2$, has
been calculated by M. Flensted-Jensen (\cite{R1}, pp. 258) to be
\begin{equation}\label{751}
e(X,Y) = \biggl( 4 \, \frac{u}{{\rm sh} \, u} \frac{v}{{\rm sh} \,
v} \frac{w}{{\rm sh} \, w}  \frac{ {\rm ch} \, w - {\rm ch} \, (u-v)
}{w^2 - (u-v)^2} \frac{ {\rm ch} \, w - {\rm ch} \, (u+v) }{w^2 -
(u+v)^2} \biggr)^{(n-3)/2}
\end{equation}
where $u, \, v$ and $w$ are the norms of $X$, $Y$ and $X+Y$,
respectively.  Note that this expression is symmetric in $X$ and
$Y$, and is only equal to 1 when $n=3$.  (\ref{751}) is also the
non-compact analogue for the expression for the (globally defined)
$e$-function given for $n$-dimensional sphere in \cite{D1} and
\cite{D2}.  We have not been able to find a proof of
Flensted-Jensen's calculation, though it should be analogous to that
in \cite{D2} by replacing $\sin$ with $\sinh$.\\

In light of our results in for compact symmetric spaces, we
conjecture that wrapping the solution to the heat equation with this
potential yields the heat kernel on the hyperbolic spaces. However,
we have not been able to calculate the heat kernel with these
potentials on $\p \cong \R^n$. Despite this, we are guaranteed the
existence of a fundamental solution whose wrap is the heat kernel on
$X$ by H\"ormander's theorem, since
$$
\frac{\partial}{\partial t} - (L_\p + \Omega_*)
$$
is hypoelliptic on $\p$.  Therefore, $\tfrac{\partial}{\partial t} -
(L_\p + \Omega_*)$ has a fundamental solution $\p$.  This
fundamental solution may be wrapped to $X$, which by Theorem
\ref{peturbHK2} is the heat kernel on $X$. $\phantom{abc} \square$

\subsection{Other non-compact semisimple Lie groups}

Extending the wrapping theorem and our results on wrapping Brownian
motion and heat kernels to other spaces may prove quite difficult.
Here, we cite the example of $SL(2,\R)$.  The (co)adjoint orbits for
$SL(2,\R)$ consist of one- and two-sheeted hyperboloids  -  the
reader is referred to \cite{HAR} for a detailed survey of $SL(2,\R)$
and the orbit method. Convolving these (co)adjoint orbits is
therefore more difficult to
describe.\\

Wrapping Brownian motion should follow in a similar way to the
procedure above.  We would then need to calculate the wrap of the
heat kernel from $\g$ to $G$.  Note for $SL(2,\R)$ the Iwasawa
decomposition $G = KAN$ is given by
$$
K =
\begin{pmatrix}
& \cos \theta & \sin \theta\\
& -\sin \theta & \cos \theta
\end{pmatrix}, \phantom{abc}
A =
\begin{pmatrix}
& a & 0\\
& 0 & a^{-1}
\end{pmatrix}, \phantom{abc}
N =
\begin{pmatrix}
& 1 & t\\
& 0 & 1
\end{pmatrix}.
$$
(see, for example, \cite{HAR}, or \cite{KNP} Ch. VI, $\S$6).  All
the elements may be conjugated into the two distinct abelian
subgroups $K$ and $A$, isomorphic to $\T$ and $\R^+$, respectively.
This creates a problem when one asks how to ``wrap" a function or
distribution.\\

In \cite{DW2}, the wrap of $j$ times an Ad-invariant Schwartz
function to a compact Lie group is calculated by summing over an
integer lattice:
\begin{equation}
\Phi (j \mu)(\exp H) = \sum_{\gamma \in \Gamma} \mu (H + \gamma),
\phantom{abcde} H \in A.
\end{equation}

We showed in section 8.3 that the wrap of $j$ times an Ad-invariant
Schwartz function to a complex Lie group may be calculated by
``bending" it:
\begin{equation}
\Phi (j \mu)(\exp H) = \mu (H), \phantom{abcde} H \in A^+.
\end{equation}

These follow from the fact that in the case of a compact Lie group,
$A \cong \T^l$, and for a complex Lie group, $A^+ \cong (\R^+)^l$.
 Since all the elements of $SL(2,\R)$ may be conjugated into the two
distinct abelian subgroups, isomorphic to $\T$ and $\R^+$, do we sum
over a lattice, or do we ``bend" to compute the wrap? Arguably, we
would need to do some combination of both.

{\sc School of Mathematics, UNSW, Kensington 2052 NSW, Australia,

And

Group Market Risk, National Australia Bank, 24/255 George St, Sydney
2000 NSW, Australia.}\\

Email: {\tt David.G.Maher@nab.com.au}


\begin{thebibliography}{12}
\small


\bibitem{A}
{\sc Anderson, A.} {\it Operator method of finding new progators
from old,}  Phys. Rev. D., {\bf 37}(2):536-539, [1987]
\bibitem{AA}
{\sc Anderson, A. and Anderson, S. B.} {\it Phase space path
integration of integrable quantum systems,}  Ann. Phys., {\bf
199}:155-186, [1990]
\bibitem{AC}
{\sc Anderson, A. and Camporesi, R.} {\it Intertwining operators for
solving differential equations, with applications to symmetric
spaces,}  Comm. Math. Phys., {\bf 130}:61-82, [1990]
\bibitem{AST}
{\sc Andler, M., Sahi, S., Torossian, C.}  {\it Convolutions of
invariant distributions:  Proof of the Kashiwara-Vergne conjecture},
Lett. Math. Phys., {\bf 69}:177-203, [2004]
\bibitem{AP}
{\sc Applebaum, D.,} {\it L\'evy processes - from probability to
finance and quantum groups,} Notices. Amer. Math. Soc. {\bf
51}(11):1336-1347, [2004]
\bibitem{ARE}
{\sc Arede, M. T.} {\it Heat kernels on Lie groups,}  Stoch. Ana.
Appl., {\bf 26}:52-62, [1991]
\bibitem{BG}
{\sc Baras, P., and Goldstein, J. A.,}  {\it The heat equation with
a singular potential,}  Trans. Amer. Math. {\bf 284}:121-139, [1966]
\bibitem{BER}
{\sc Berger, M.}  {\it Geometry of the spectrum I,} Proc. Sympos.
Pure Math. {\bf 27},2, AMS, Providence, R.I. pp. 265-283, [1975]
\bibitem{BGM}
{\sc Berger, M., Gauduchon, P., and Mazet, E.}  {\it Le spectre
d'une vari\'et\'e Riemannienne,} Lecture notes in Math. {\bf 194}
Springer-Verlag, Berlin, [1971]
\bibitem{CAM}
{\sc Camporesi, R.} {\it Harmonic analysis and propagators on
homogeneous spaces,}  Phys. Rep., {\bf 196}:1-134, [1990]
\bibitem{CHU}
{\sc Chung, K.H.}  {\it Compact group actions and harmonic
analysis,}  PhD Thesis, University of New South Wales, [1999]
\bibitem{D1}
{\sc Dooley, A.H.}  {\it Orbital convolutions, wrapping maps and
$e$-functions,}  Proc. CMA, ANU, [2002]
\bibitem{D2}
{\sc Dooley, A.H.}  {\it Global versions of the $e$-function for
compact symmetric spaces,}  In preparation.
\bibitem{DW1}
{\sc Dooley, A.H. and Wildberger, N.J.}  {\it Global character
formulae for compact Lie groups}  Trans. Amer. Math. Soc., {\bf
351}(2):477-495, [1999]
\bibitem{DW2}
{\sc Dooley, A.H. and Wildberger, N.J.}  {\it Harmonic Analysis and
the Global Exponential Map for Compact Lie Groups,} Funktsional.
Anal. i Prilozhen. \textbf{27}(1):25-32; Eng. Trans., [1993].
Funct. Ana. Appl. \textbf{27}:21-27; MR \textbf{94e:}22032, [1993]
\bibitem{DRW}
{\sc Dooley, A.H., Repka, J., and Wildberger, N.J.}  {\it Sums of
adjoint orbits,}  Lin. Multilin. Alg., {\bf 36}:79-101, [1993]
\bibitem{DOW}
{\sc Dowker, J. S.} {\it When is the `sum over classical paths'
exact?}  J. Phys. A, {\bf 3}:451-461, [1970]
\bibitem{DUF}
{\sc Duflo, M.} {\it Op\'erateurs diff\'erential bi-invariants sur
un groupe de Lie,} Ann. Sci. \'Ecole Norm. Sup., {\bf 10}:265-288,
[1977]
\bibitem{GAN}
{\sc Gangolli, R.}  {\it Asymptotic behaiviour of spectra of compact
quotients of certain symmetric spaces}, Acta. Math.
\textbf{121}:151-192, [1968]
\bibitem{GAN2}
{\sc Gangolli, R.}  {\it Spherical Functions on a Semisimple Lie
Group}, In: ``Symmetric Spaces: Short courses presented at
Washington University", pp. 41-92, Marcel Decker, Inc., New York,
[1972]
\bibitem{GV}
{\sc Gangolli, R. and Varadarajan, V.S.,}  {\it Harmonic Analysis of
Spherical Functions on Real Reductive Groups}, Springer-Verlag, New
York, [1988]
\bibitem{HS}
{\sc Hall, B. C. and Stenzel, M. B.}  {\it Sharp bounds for the heat
kernel on certain symmetric spaces of non-compact type}  Contemp.
Math., {\bf 317}:117-135, [2003]
\bibitem{HAR}
{\sc Harinck, P.}  {\it Orbit method for $SL(2,\R)$}   European
women in mathematics (Trieste, 1997), Hindawi Publ. Corp., Stony
Brook, NY, pp 113-121, [1999]
\bibitem{HE1}
{\sc Helgason, S.}  {\it Differential Geometry, Lie Groups and
Symmetric Spaces,}  Academic Press, [1978]
\bibitem{HE2}
{\sc Helgason, S.}  {\it Groups and Geometric Analysis,}
Mathematical Surveys and Monographs, {\bf 83}, AMS, [2000]
\bibitem{HE3}
{\sc Helgason, S.}  {\it Geometric analysis on symmetric spaces,}
Mathematical Surveys and Monographs, {\bf 39}, AMS, [1994]
\bibitem{KV}
{\sc Kashiwara, M., Vergne, M.} {\it The Campbell-Hausdorff formula
and invariant hyperfunctions,}  Invent. Math. {\bf 47}:249-272,
[1978]
\bibitem{KIN}
{\sc Kingman, J.F.C.}  {\it Random Walks with Spherical Symmetry},
Acta Math., {\bf 109}:11-53, [1963]
\bibitem{KNP}
{\sc Knapp, A. W.}  {\it Lie Groups, Beyond an Introduction,}
Birkh\"auser, Boston, [2002]
\bibitem{KON}
{\sc Kontsevich, M.} {\it Deformation quantization of Poisson
manifolds,}  Lett. Math. Phys. {\bf 66}(3):157-216, [2003]
\bibitem{LIAO}
{\sc Liao, M.}  {\it L\'evy processes and Fourier analysis on
compact Lie groups,}  Ann. of Prob.  {\bf 32}(2):1553-1573, [2004]
\bibitem{LIAO2}
{\sc Liao, M.}  {\it L\'evy processes in Lie groups,} Cambridge
University Press, [2004]
\bibitem{LOW}
{\sc Low, S.}  {\it Path integration on space-times with symmetry,}
Ph.D thesis, University of Texas at Austin, [1985]
\bibitem{M}
{\sc Maher, D. G.,}  {\it Brownian motion and heat kernels on
compact Lie groups and symmetric spaces,}  Ph.D thesis, University
of New South Wales.
\bibitem{M2}
{\sc Maher, D. G.,}  {\it Wrapping Brownian motion and heat kernels
I: compact Lie groups,}  Preprint.
\bibitem{M3}
{\sc Maher, D. G.,}  {\it The wrapping map, wave equation, and
Huygens principal,}  In preparation.
\bibitem{OP}
{\sc Olshanetski, M. A., Perelomov, A.M.}  {\it  Quantum integrable
systems related to Lie algebras,} Phys. Rep., {\bf 94}:312-404,
[1993]
\bibitem{PI}
{\sc Peak, D. and Inomata, A.} {\it Summation over Feynman histories
in polar coordinates,}  J. Math. Phys., {\bf 10}(9):1422-1428,
[1969]
\bibitem{ROU}
{\sc Rouvi\`ere, F.} {\it D\'emonstration de la conjecture de
Kashiwara-Vergne pour $SL_2(\R)$}, C. R. Acad. Sci. Paris S\'erie I,
{\bf 292}657660, [1981]
\bibitem{R1}
{\sc Rouvi\`ere, F.} {\it Invariant analysis and contractions of
symmetric spaces, Part I,}  Compositio Math. {\bf 73}:241-270,
[1990]
\bibitem{R2}
{\sc Rouvi\`ere, F.} {\it Invariant analysis and contractions of
symmetric spaces, Part II,}  Compositio Math. {\bf 80}:111-136,
[1991]
\bibitem{T}
{\sc Torossian, C.}  {\it M\'ethodes de Kashiwara-Vergne-Rouvi\`ere
pour les espaces sym\'etriques,}  Noncommutative harmonic analysis,
pp 459-486, Progr. Math., 220, Birkh\"auser, Boston, [2004]
\bibitem{WAT}
{\sc Watanabe, S.}  {\it L\'evy's stochastic area formula and
Brownian motion on compact Lie groups,}  In: ``It\^o's Stochastic
Calculus and Probability Theory", pp 401-411.  Springer-Verlag,
[1996]
\bibitem{WLD}
{\sc Wildberger, N.J.}  {\it Hypergroups, Symmetric Spaces, and
Wrapping Maps},  Probability Measures on Groups and Related
Structures, vol XI of Proceedings Oberwofach,  World Scientific
Publishing Co. Pty. Ltd, Singapore, [1994]





\end{thebibliography}
\end{document}